\documentclass[12pt]{amsart}
\usepackage{graphicx}
\usepackage{rotate,epic,eepic,epsf,epsfig}

\usepackage{amsmath}
\usepackage{amsfonts}
\usepackage{amssymb}
\newtheorem{defi}{D\'{e}finition}[section]
\newtheorem{prop}{Proposition}[section]
\newtheorem{theo}{Th\'eor\`eme}[section]
\newtheorem{lem}{Lemme}[section]

\newtheorem{exe}{Exemple}[section]

\newcommand{\resumename}{R\'esum\'e}
\newenvironment{resume}{\narrower\footnotesize\bf
\noindent\resumename.\quad\footnotesize\rm}{\par\bigskip}

\DeclareFixedFont{\petitefonte}{\encodingdefault}%
{\familydefault}{\seriesdefault}{\shapedefault}{5pt}

\begin{document}


\title[C\^one diamant symplectique]{Le c\^one diamant symplectique}
\author[D. Arnal et O. Khlifi]{Didier Arnal et Olfa Khlifi}
\address{
Institut de Math\'ematiques de Bourgogne\\
UMR CNRS 5584\\
Universit\'e de Bourgogne\\
U.F.R. Sciences et Techniques
B.P. 47870\\
F-21078 Dijon Cedex\\France} \email{Didier.Arnal@u-bourgogne.fr}
\address{D\'epartement de Mathematiques, Facult\'e des Sciences de Sfax, Route de Soukra, km 3,5,
B.P. 1171, 3000 Sfax, Tunisie.} \email{khlifi\_olfa@yahoo.fr}

\begin{abstract}
The diamond cone is a combinatorial description for a basis in a indecomposable module for the nilpotent factor $\mathfrak n^+$ of a semi simple Lie algebra. After N.J. Wildberger who introduced this notion for $\mathfrak{sl}(3)$, this description was achevied in \cite{ABW} for $\mathfrak{sl}(n)$ and in \cite{AAK} for the rank 2 semi-simple Lie algebras.\\

In the present work, we generalize these constructions to the Lie algebras $\mathfrak{sp}(2n)$. The symplectic semi-standard Young tableaux were defined by C. de Concini in \cite{DeC}, they form a basis for the shape algebra of $\mathfrak{sp}(2n)$.
We introduce here the notion of symplectic quasi-standard Young tableaux, these tableaux give the diamond cone for $\mathfrak{sp}(2n)$.\\
\end{abstract}




\thanks{
Ce travail a \'et\'e effectu\'e dans le cadre de l'accord CMCU 06 S 1502, O. Khlifi remercie l'Universit\'e de Bourgogne pour l'accueil dont elle a b\'en\'efici\'e au cours de ses s\'ejours, D. Arnal remercie la Facult\'e des Sciences de Sfax pour l'accueil dont il a b\'en\'efici\'e au cours de ses s\'ejours.}


\maketitle \vspace{.20cm}

\begin{resume}
Si $\mathfrak n^+$ est le facteur nilpotent d'une alg\`ebre semi-simple $\mathfrak g$, le c\^one diamant de $\mathfrak g$ est la description combinatoire d'une base d'un $\mathfrak n^+$ module ind\'ecomposable naturel. Cette notion a \'et\'e introduite par N. J. Wildberger pour $\mathfrak{sl}(3)$, le c\^one diamant de $\mathfrak{sl}(n)$ est d\'{e}crit dans \cite{ABW}, celui des alg\`ebres semi simples de rang 2 dans \cite{AAK}.\\

Dans cet article, nous g\'en\'eralisons ces constructions au cas des alg\`ebres de Lie $\mathfrak{sp}(2n)$. Les tableaux de Young semi standards symplectiques ont \'et\'e d\'efinis par C. de Concini dans \cite{DeC}, ils forment une base de l'alg\`ebre de forme de $\mathfrak{sp}(2n)$. Nous introduisons ici la notion de tableaux de Young quasi standards symplectiques, ces derniers d\'ecrivent le c\^one diamant de $\mathfrak{sp}(2n)$.\\
\end{resume}

\medskip
\section{Introduction}

\

Soit $\mathfrak g$ une alg\`ebre de Lie semi simple complexe de dimension finie et
$$
\mathfrak g=\mathfrak h+\sum_{\alpha \in \Phi}\mathfrak
g^\alpha=\mathfrak h+\sum_{\alpha \in \Phi}\mathbb C X_\alpha
$$
sa d\'ecomposition en sous-espaces radiciels.\\

La th\'eorie des modules simples de dimension finie de $\mathfrak g$ est tr\`es bien connue et assez explicite. Ayant fix\'e un syst\`eme de racines simples $\Delta$, on note $\Phi^+$ l'ensemble des racines positives, $\mathfrak n^+=\sum_{\alpha\in\Phi^+}\mathfrak g^\alpha$, on sait qu'un tel module simple $V^\lambda$ est caract\'eris\'e \`a \'equivalence pr\`es par son plus haut poids $\lambda$ qui est entier et dominant. Soit $\Lambda$ l'ensemble des poids entiers dominants. Cette th\'eorie peut se r\'esumer \`a la description de l'alg\`ebre de forme de $\mathfrak g$. Cette alg\`ebre est l'espace
$$
\mathbb V=\bigoplus_{\lambda \in \Lambda}V^\lambda,
$$
muni d'une multiplication associative et commutative naturelle (voir \cite{FH}).\\

Un probl\`eme combinatoire classique est alors de d\'ecrire explicitement cette alg\`ebre, en particulier d'en donner une base, form\'ee d'une union de bases de chaque $V^\lambda$. Par exemple dans le cas o\`u $\mathfrak g=\mathfrak{sl}(n)$, on note $\mathbb S^\lambda$ le module $V^\lambda$ et $\mathbb S^\bullet$ l'alg\`ebre de forme de $\mathfrak{sl}(n)$. Cette alg\`ebre s'identifie \`a une alg\`ebre de fonctions polyn\^omes sur le groupe de Lie $SL(n)$, et on conna\^\i t depuis le $19^{\grave eme}$ si\`ecle une telle base. On peut l'indexer par l'ensemble $SS^\bullet$ des tableaux de Young semi standards, remplis par des coefficients dans $\{1,\dots,n\}$ et dont les colonnes sont de hauteur inf\'erieure \`a $n$ : chaque tableau d\'efinit naturellement une fonction polyn\^ome, produits de sous d\'eterminants sur $SL(n)$, ces fonctions forment une base de $\mathbb S^\bullet$, les tableaux de forme $\lambda$ d\'efinissant une base de $V^\lambda$.\\

Dans la suite, on notera $v_{-\lambda}$ un vecteur de plus bas poids de $V^\lambda$, $V^\lambda$ est engendr\'e par l'action de $\mathfrak n^+$ sur le vecteur $v_{-\lambda}$. On notera $V^\lambda_{\mathfrak n^+}$ l'espace $V^\lambda$ vu comme un $\mathfrak n^+$ module monog\`ene. Ces modules monog\`enes sont en fait maximaux et caract\'eris\'es par les nombres entiers naturels $a_\alpha$ ($\alpha \in \Delta$) tels que
$$
X_\alpha^{a_\alpha}v_{-\lambda}\neq0\quad\text{et}\quad X_\alpha^{a_\alpha+1}v_{-\lambda}=0\quad(\alpha \in \Delta).
$$
La description des modules monog\`enes nilpotents de $\mathfrak n^+$ semble donc se r\'esumer \`a la description d'une nouvelle alg\`ebre $\mathbb V_{red}$, quotient de l'alg\`{e}bre de forme $\mathbb V$ et que l'on appellera l'alg\`ebre de forme r\'eduite de $\mathfrak g$. Cette alg\`ebre ne sera plus la somme directe des $V^\lambda_{\mathfrak n^+}$ mais en fait un $\mathfrak n^+$ module ind\'ecomposable, union de tous ces modules, avec la stratification naturelle : $V^\mu_{\mathfrak n^+}\subset V^\lambda_{\mathfrak n^+}$ si et seulement si $\mu\leq\lambda$.\\

Le probl\`eme combinatoire est maintenant de d\'ecrire une base de l'alg\`ebre de forme r\'eduite, adapt\'ee \`a la stratification, c'est \`a dire une base union de bases des $V^\lambda_{\mathfrak n^+}$, la base de $V^\lambda_{\mathfrak n^+}$ contenant toutes celles des $V^\mu_{\mathfrak n^+}$ si $\mu\leq \lambda$.\\

Supposons de nouveau que $\mathfrak g=\mathfrak{sl}(n)$. Ce cas a \'et\'e \'etudi\'e dans \cite{ABW}. L'alg\`ebre de forme r\'eduite, $\mathbb S^\bullet_{red}$ est isomorphe \`a l'alg\`ebre $\mathbb C[N^+]$ des fonctions polyn\^omes sur le groupe $N^+=\exp \mathfrak n^+$. Une base adapt\'ee de cette alg\`ebre est donn\'ee par l'ensemble $QS^\bullet$ des tableaux de Young appel\'es quasi standards. La base de $V^\lambda_{\mathfrak n^+}$ \'etant donn\'ee par l'ensemble des tableaux de Young quasi standards de forme inf\'erieure ou \'egale \`a $\lambda$. Autrement dit, l'ensemble des tableaux de Young quasi standards de forme $\lambda$ forme une base d'un suppl\'ementaire de $\sum_{\mu<\lambda}V^\mu_{\mathfrak n^+}$ dans $V^\lambda_{\mathfrak n^+}$. En reprenant la terminologie de N. J. Wildberger, on dit qu'on a d\'ecrit le c\^one diamant de $\mathfrak{sl}(n)$, \cite{W}.\\

R\'ecemment, avec B. Agrebaoui, nous avons r\'ealis\'e la m\^eme construction combinatoire pour les alg\`ebres de rang 2 : $\mathfrak{sl}(2)\times\mathfrak{sl}(2)$, $\mathfrak{sl}(3)$, $\mathfrak{sp}(4)$ et $g_2$. Nous avons ainsi d\'ecrit leur c\^one diamant en utilisant pour chacune d'elles la bonne notion de tableau de Young quasi standard (\cite{AAK}).\\

Le but de cet article est de traiter de la m\^eme fa\c con les alg\`ebres de Lie symplectiqes $\mathfrak{sp}(2n)$. Pour ces alg\`ebres, on peut d\'efinir la classe des tableaux de Young semi standards, ce qui donne une base de l'alg\`ebre de forme $\mathbb S^{\langle\bullet\rangle}$ de $\mathfrak{sp}(2n)$.\\

On regarde d'abord $\mathfrak{sp}(2n)$ comme une sous alg\`ebre de Lie de l'alg\`ebre $\mathfrak{sl}(2n)$, de telle fa\c con que, avec des notations \'evidentes,
$$
\mathfrak h_{\mathfrak{sp}(2n)}=\mathfrak h_{\mathfrak{sl}(2n)}\cap \mathfrak{sp}(2n),\quad \mathfrak n^+_{\mathfrak{sp}(2n)}=\mathfrak n^+_{\mathfrak{sl}(2n)}\cap \mathfrak{sp}(2n).
$$
On se limite alors aux $\mathfrak{sl}(2n)$ modules simples $\mathbb S^\lambda$ qui correspondent aux tableaux de Young n'ayant pas de colonne de hauteur $>n$, la res\-triction de $\mathbb S^\lambda$ \`a $\mathfrak{sp}(2n)$ contient exactement un $\mathfrak{sp}(2n)$ module simple $\mathbb S^{\langle \lambda \rangle}$ de plus haut poids $\lambda|_{\mathfrak h_{\mathfrak{sp}(2n)}}$. On d\'ecrit ainsi exactement l'ensemble des $\mathfrak{sp}(2n)$ modules simples et il existe deux combinatoires, celles de de Concini et celle de Kashiwara-Nakashima qui permettent de s\'election\-ner, parmi les tableaux de Young semi standards de forme $\lambda$, une base de $V^{\langle \lambda\rangle}$. On dira que ces tableaux sont les tableaux semi standards symplectiques, voir \cite{DeC}, \cite{KN}.\\

Dans cet article, on va d\'efinir la notion de tableau de Young quasi standard symplectique et montrer que, comme dans le cas de $\mathfrak{sl}(n)$ ou des alg\`ebres de rang 2, les tableaux quasi standards symplectiques de forme $\lambda$ forment une base d'un suppl\'ementaire de $\sum_{\mu<\lambda}\mathbb S^{\langle\mu\rangle}_{\mathfrak n^+_{\mathfrak{sp}(2n)}}$ dans $\mathbb S^{\langle\lambda\rangle}_{\mathfrak n^+_{\mathfrak{sp}(2n)}}$. On obtiendra ainsi une base de l'alg\`ebre de forme r\'eduite de $\mathfrak{sp}(2n)$. Cette alg\`ebre, not\'ee $\mathbb S^{\langle\bullet\rangle}_{red}$, est isomorphe \`a $\mathbb C[N^+_{\mathfrak{sp}(2n)}]$ et sa structure de $\mathfrak n^+_{\mathfrak{sp}(2n)}$ module ind\'ecomposable est bien d\'ecrite par notre base. On aura ainsi d\'ecrit le c\^one diamant des alg\`ebres $\mathfrak{sp}(2n)$.\\

\section{Tableaux de Young semi et quasi standards pour $\mathfrak{sl}(n)$}

\

Dans cette section, on va rappeler les d\'{e}finitions, les notations et les r\'{e}sultats de l'article \cite{ABW} qui \'etudie le cas des alg\`ebres $\mathfrak{sl}(n)$. On esquissera aussi une nouvelle preuve du r\'esultat principal de ce travail, en utilisant le jeu de taquin de Sch\"utzenberger. C'est cette preuve qui sera g\'en\'eralis\'ee pour $\mathfrak{sp}(2n)$.\\

\subsection{Tableaux de Young semi standards pour $\mathfrak{sl}(n)$}

\

L'alg\`{e}bre de Lie $\mathfrak{sl}(n)$ est l'ensemble des matrices complexes carr\'{e}es d'ordre $n$ et de trace nulle.
Le groupe de Lie correspondant, $SL(n)$, est l'ensemble des matrices carr\'{e}es d'ordre $n$ et de d\'{e}terminant 1.\\

L'ensemble $\mathfrak h$ des matrices diagonales $H=diag(\kappa_1,\dots,\kappa_n)$ (avec $\sum_i \kappa_i=0$) est une sous alg\`{e}bre de Cartan de $\mathfrak{sl}(n)$. On d\'{e}finit les formes lin\'{e}aires $\theta_i$ sur $\mathfrak h$ en posant $\theta_i(H)=\kappa_i$. On choisit l'ensemble des racines simples $\Delta=\{\alpha_i=\theta_{i+1}-\theta_i,~~1\leq i<n\}$. Pour $1\leq k<n$, l'action naturelle de $\mathfrak{sl}(n)$ sur $\wedge^k\mathbb C^n$ d\'efinit des modules irr\'eductibles de plus haut poids $\omega_k=\theta_1+\dots+\theta_k$. Ces modules sont les repr\'{e}sentations fondamentales de $\mathfrak{sl}(n)$.\\

Chaque $\mathfrak{sl}(n)$ module simple est caract\'{e}ris\'{e} par son plus haut poids
$$
\lambda=\sum _{k=1} ^{n-1} a_k\omega_k
$$
o\`u les $a_k$ sont des entiers naturels. Notons ce module irr\'{e}ductible $\mathbb S^\lambda$, c'est un sous module de
$$
Sym ^{a_1}(\mathbb C^n)\otimes  Sym ^{a_2} (\wedge ^2 \mathbb C^n)\otimes \dots \otimes Sym ^{a_{n-1}} (\wedge ^{n-1}\mathbb{C}^n).
$$
La th\'eorie classique des $\mathfrak{sl}(n)$ modules simples dit que l'ensemble des modules simples est en bijection avec l'ensemble $\Lambda$ des poids entiers positifs, et que l'application $\lambda\mapsto (a_1,\dots,a_n)$ est une bijection de $\Lambda$ sur $\mathbb N^{n-1}$.\\

Soit $(e_1,\dots,e_n)$ la base canonique de $\mathbb C^n$. Le d\'{e}terminant de la sous
matrice de $g$ obtenue en ne consid\'{e}rant que les lignes
$i_1,\dots,i_k$ et les colonnes $j_1, \dots,j_d$ est not\'e $\det~(g;i_1,\dots,i_k;j_1,\dots,j_k)$. Une base de
$\mathbb S^{\omega_k}$ est donn\'{e}e par l'ensemble des fonctions
sous d\'eterminant suivantes~:
$$\begin{aligned}
\delta^{(k)}_{i_1,\dots,i_k}(g)&=\det~(g;i_1,\dots,i_k;1,\dots,k)\\
&=\langle e^\star_{i_1}\wedge \dots \wedge e^\star_{i_k}, ge_1\wedge \dots \wedge ge_k\rangle
\end{aligned}
$$
o\`{u} $g \in SL(n)$, et $i_1< i_2<\dots<i_k$.

On note cette fonction par une colonne :
$$
\delta^{(k)}_{i_1,\dots,i_k}=\begin{array}{|c|} \hline
i_1 \\
\hline
i_2 \\
\hline
\vdots \\
\hline
i_k \\
\hline
\end{array}.
$$
Si $i_1=1,i_2=2,\dots,i_k=k$, la colonne sera dite triviale.

Le groupe $SL(n)$ agit sur ces colonnes par l'action r\'eguli\`ere gauche :
$$
(g.\delta^{(k)}_{i_1,\dots,i_k})(g')=\delta^{(k)}_{i_1,\dots,i_k}(^tgg').
$$
Par construction, cette action co\"\i ncide avec l'action naturelle de $SL(n)$ sur $\wedge^k\mathbb C^n$. La colonne triviale est le vecteur de poids $\omega_k$, on la choisit comme le vecteur de plus haut poids de $\mathbb S^{\omega_k}$, ce module est maintenant d\'efini univoquement (pas \`a un op\'erateur scalaire pr\`es).

On notera un produit de fonctions $\delta$ comme un tableau, form\'e d'une juxtaposition de colonnes qu'on appellera tableau de Young. Un tableau de Young vide $T$ est une suite finie de colonnes $c_1,\dots,c_r$. Chaque colonne verticale $c_j$ est form\'ee de $\ell_j$ cases vides. Ces cases sont rep\'er\'ees par un double indice : pour la colonne $c_j$, ce sont les cases $(1,j)$,\dots, $(\ell_j,j)$. On suppose $1\leq \ell_r\leq \dots\leq \ell_1\leq n-1$. La forme du tableau $form(T)$ est le $n-1$ uplet $(a_1,\dots,a_{n-1})$ s'il y a $a_1$ colonnes de hauteur 1,\dots, $a_{n-1}$ colonnes de hauteur $n-1$. On remplit le tableau avec des entiers $t_{ij}$ plac\'es dans les cases vides.\\

Ainsi, l'ensemble des tableaux de Young forme une base de l'alg\`{e}bre sym\'etrique :
$$\aligned
Sym^{\bullet}(\bigwedge\mathbb C^n)&=Sym^{\bullet} (\mathbb{C}^{n} \oplus \wedge^{2} \mathbb{C}^{n}\oplus \dots \oplus \wedge^{n-1} \mathbb{C}^{n})\\
&=\sum_{a_1,\dots,a_{n-1}}Sym ^{a_1}(\mathbb C^n)\otimes \dots \otimes Sym ^{a_{n-1}} (\wedge ^{n-1}\mathbb{C}^n).
\endaligned
$$
Si $\lambda=\sum a_k\omega_k$, le module $\mathbb S^\lambda$ est alors \'equivalent au sous-module de $Sym^{\bullet}(\bigwedge\mathbb C^n)$ engendr\'e par l'action de $\mathfrak{sl}(n)$ sur le tableau de Young $T^\lambda$ ayant exactement $a_1$ colonnes triviales de hauteur 1, \dots , $a_{n-1}$ colonnes triviales de hauteur $n-1$.\\

Soit $N^+$ le groupe des matrices $n\times n$ triangulaires sup\'erieures avec des 1 sur la diagonales. On montre que l'alg\`ebre des fonctions polynomiales en les coefficients de $g\in SL(n)$ $N^+$ invariantes par multiplication \`a droite est engendr\'ee par les fonctions $\delta^{(k)}_{i_1,\dots,i_k}$. Cette alg\`ebre est donc un quotient de l'alg\`ebre $Sym^\bullet(\bigwedge \mathbb C^n)$. En tant que $\mathfrak{sl}(n)$ module, elle est engendr\'ee par les fonctions $T^\lambda$, c'est la somme directe des $\mathbb S^\lambda$.\\

\begin{defi}

\

L'alg\`{e}bre de forme de $SL(n)$ est le $\mathfrak{sl}(n)$ module :
$$
{\mathbb S}^\bullet=\bigoplus_{\lambda\in \Lambda}~\mathbb S^\lambda
$$
vu comme le quotient de $Sym^\bullet(\bigwedge \mathbb C^n)$ d\'efini ci-dessus.

Un tableau de Young de forme $\lambda=(a_1,\dots,a_{n-1})$ est dit semi standard si son remplissage se fait par des entiers $\leq n$ qui sont croissants de gauche \`{a} droite le long de chaque ligne et strictement croissants de haut en bas le long de chaque colonne.\\
\end{defi}

La th\'eorie classique des tableaux de Young semi standards dit que ces tableaux forment une base de l'espace $\mathbb S^\bullet$. Plus pr\'ecis\'ement :\\
\begin{theo}

\

\begin{itemize}
\item[1)] On a les isomorphismes d'alg\`ebre :
$$
{\mathbb S}^\bullet\simeq \mathbb{C}[SL(n)]^{N^+} \simeq \mathbb{C}[\delta^{(k)}_{i_1,\dots,i_k}]/_{\hbox{$\mathcal{PL}$ }}.$$
L'id\'{e}al $\mathcal{PL}$ est l'id\'eal engendr\'e par les relations Pl\"ucker : pour $p \geq q \geq r$,
$$\aligned
&0=\delta^{(p)}_{i_1,i_2,\dots,i_p}\delta^{(q)}_{j_1,j_2,\dots,j_q}+\\
&\hskip 1cm+\sum\limits_{\begin{array}{c}
A \subset\{i_1,\dots,i_p\}\\
\# A = r
\end{array}}\pm \delta^{(p)}_{(\{i_1,\dots,i_p\}\setminus A)\cup\{j_1,\dots,j_r\}}\delta^{(q)}_{A\cup \{j_{r+1},\dots,j_q\}}.\endaligned
$$
\item[2)] Si $\lambda=a_1\omega_1+\dots+a_{n-1}\omega_{n-1}=(a_1,\dots,a_{n-1})$, alors une base de ${\mathbb S}^\lambda$ est donn\'{e}e par l'ensemble des tableaux de Young semi standards de forme $\lambda$.
\item[3)] La relation d'ordre sur les poids $\mu\leq\lambda$ corespond \`a la relation d'ordre partielle $b_k\leq a_k$ pour tout $k$ si $\mu=(b_1,\dots,b_k)$ et $\lambda=(a_1,\dots,a_k)$.\\
\end{itemize}
\end{theo}

\begin{exe}

Pour le cas de $\mathfrak{sl}(3)$ ($n=3$), on a une seule relation
de Pl\"ucker:
$$
\renewcommand{\arraystretch}{0.7}{\begin{array}{l}
\framebox{$1$}\framebox{$3$}\\
\framebox{$2$} \\
\end{array}}+ \renewcommand{\arraystretch}{0.7}{\begin{array}{l}
\framebox{$2$}\framebox{$1$}\\
\framebox{$3$} \\
\end{array}} - \renewcommand{\arraystretch}{0.7}{\begin{array}{l}
\framebox{$1$}\framebox{$2$}\\
\framebox{$3$} \\
\end{array}}=0.
$$
L'alg\`ebre de forme ${\mathbb S}^\bullet$, lorsque $n=3$, est une
sous alg\`ebre de $Sym^\bullet(\bigwedge \mathbb C^3)$, on a vu comment
d\'{e}finir une base de cette derni\'{e}re, form\'{e}e de tableaux
de Young. La base de $\mathbb S^\bullet$ est obtenue en
\'{e}liminant les tableaux non semi standards. C'est \`{a} dire
exactement ceux qui contiennent le sous tableau
$\renewcommand{\arraystretch}{0.7}{\begin{array}{l}
\framebox{$2$}\framebox{$1$}\\
\framebox{$3$} \\
\end{array}}$.\\
\end{exe}

\

\subsection{Tableaux de Young Quasi-standards pour $\mathfrak{sl}(n)$}

\

Pour construire l'alg\`ebre de forme r\'eduite \`a partir de l'alg\`ebre de forme, on restreint les fonctions polynomiales $N^+$ invariantes sur $SL(n)$ au sous groupe $N^-=~^tN^+$.\\

\begin{defi}

\

On appelle alg\`ebre forme r\'eduite, et on note $\mathbb{S}^\bullet_{red}$, le quotient :
$$
\mathbb{S}^\bullet_{red}=\mathbb{S}^\bullet/ \big<\delta^{(k)}_{1,\dots, k}-1~~ \big >.
$$
\end{defi}

\begin{theo}

\

En tant qu'alg\`ebre, $\mathbb{S}^\bullet_{red}$ est l'alg\`ebre des fonctions polynomiales sur le groupe $N^-$. C'est aussi le quotient de l'alg\`ebre sym\'etrique sur les fonctions $\delta^{(k)}_{i_1,\dots,i_k}$ non triviales ($i_k>k$) par l'id\'eal des relations de Pl\"ucker r\'eduites, c'est \`a dire des relations de Pl\"ucker dans lesquelles on suprime les colonnes triviales.

En tant que $\mathfrak n^+$ module, $\mathbb{S}^\bullet_{red}$ est ind\'ecomposable et c'est l'union des modules $V^\lambda_{\mathfrak n^+}=\mathbb{S}^\lambda_{\mathfrak n^+}$, stratifi\'ee par :
$$
\mu\leq \lambda \Longleftrightarrow \mathbb{S}^\mu_{\mathfrak n^+} \subset \mathbb{S}^\lambda_{\mathfrak n^+}.
$$

\end{theo}

\begin{defi}

\

On consid\'ere un tableau semi standard $T=(t_{ij})$. Si le haut de la premi\`ere colonne de $T$ (les $s$ premi\`eres lignes) est trivial, si $T$ contient une colonne de hauteur $s$ et si pour tout $j$ pour lequel ces entr\'ees existent, on a $t_{s(j+1)}<t_{(s+1)j}$, on dit que $T$ n'est pas quasi standard en $s$. S'il n'existe aucun tel $s$, on dit que $T$ est quasi standard.
\end{defi}

\begin{exe}

\

La relation de  Pl\"ucker r\'{e}duite pour $\mathfrak{sl}(3)$ est :
$$
\renewcommand{\arraystretch}{0.7}{
\framebox{$3$}} +
\renewcommand{\arraystretch}{0.7}{\begin{array}{l}
\framebox{$2$}\\
\framebox{$3$} \\
\end{array}}-\renewcommand{\arraystretch}{0.7}{\begin{array}{l}
\framebox{$1$}\framebox{$2$}\\
\framebox{$3$} \\
\end{array}}=0.
$$

Cette relation contient un seul tableau non quasi standard : le dernier.\\

\end{exe}

Notons $SS^\lambda$ (resp. $QS^\lambda$) l'ensemble des tableaux de Young semi standards (resp. quasi standards) de forme $\lambda$.\\

Le r\'esultat principal de \cite{ABW} est que les tableaux quasi standards d\'ecrivent le c\^one diamant de $\mathfrak{sl}(n)$. Donnons une preuve de ce r\'esultat utilisant le jeu de taquin de Sch\"utzenberger.

Soient $S$ et $T$ deux tableaux de Young vides de forme $\mu=form(S)=(b_1,\dots,b_n)\leq \lambda =form(T)=(a_1,\dots,a_{n-1})$. On place $S$ dans le coin en haut \`a gauche de $T$. Un coin int\'erieur de $S$ est une case $(x,y)$ de $S$ telle que, imm\'ediatement \`a droite et imm\'ediatement en dessous de cette case, il n'y a pas de case de $S$. Un coin ext\'erieur de $T$ est une case vide $(x',y')$ qu'on peut ajouter \`a $T$ de telle fa\c con que $T\cup\{(x',y')\}$ soit encore un tableau de Young (ses colonnes sont de hauteurs d\'ecroissantes et commencent \`a la premi\`ere ligne).

On laisse le tableau $S$ vide et on remplit le `tableau tordu' $T\setminus S$ de forme $\lambda\setminus\mu$ par des entiers $t_{ij}\leq n$ de fa\c con semi standard : pour tout $i$ et tout $j$, $t_{ij}<t_{(i+1)j}$ et $t_{ij}\leq t_{i(j+1)}$, si les cases correspondentes sont dans $T\setminus S$. On choisit un coin int\'erieur de $S$ et on l'identifie par une \'etoile : $\boxed{\star}$. On dira qu'on a un tableau tordu $T\setminus S$ point\'e. Par exemple,
$$
\begin{tabular}{|c|c|c|} \hline &$2$& \multicolumn{1}{|c|}{$4$}
\\  \hline
   $\star$   & $3$&  \multicolumn{1}{|c|}{$5$ } \\ \cline{1-3}
 $4$  & $6$  \\ \cline{1-2}
 $5$&$7$\\ \cline{1-2}
 \end{tabular}
$$
est un tableau tordu point\'e.

Le jeu de taquin consiste \`a d\'eplacer cette case $\boxed{\star}$ dans $T$. Apr\`es un certain nombre de d\'eplacements, le tableau $T$ est devenu un tableau $T'$ dans lequel la case point\'ee est \`a la place $(i,j)$. Alors
\begin{itemize}
\item[] Si la case $(i,j+1)$ existe et si la case $(i+1,j)$ n'existe pas ou $t_{(i+1)j}>t_{i(j+1)}$, on pousse $\boxed{\star}$ vers la droite, c'est \`a dire, on remplace $T'$ par le tableau $T"$ o\`u en $(i,j)$, on met $\boxed{t_{i(j+1)}}$, on met $\boxed{\star}$ en $(i,j+1)$, on ne modifie pas les autres entr\'ees de $T'$.\\

\item[] Si la case $(i+1,j)$ existe et si la case $(i,j+1)$ n'existe pas ou $t_{(i+1)j}\leq t_{i(j+1)}$, on pousse $\boxed{\star}$ vers le bas, c'est \`a dire, on remplace $T'$ par le tableau $T"$ o\`u en $(i,j)$, on met $\boxed{t_{(i+1)j}}$, on met $\boxed{\star}$ en $(i+1,j)$, on ne modifie pas les autres entr\'ees de $T'$.\\

\item[] Si les cases $(i+1,j)$ et $(i,j+1)$ n'existent pas, on supprime la case $\boxed{\star}$. La case $(i,j)$ n'est plus une case de $T"$ mais le tableau form\'e des cases de $T"$ et de la case $(i,j)$ est un tableau de Young. La case $(i,j)$ est un coin ext\'erieur de $T"$.\\
\end{itemize}

\begin{exe}

\

$$\aligned
T~&=~{\begin{tabular}{|c|c|c|}
\hline
 & $2$ & \multicolumn{1}{|c|}{$4$}\\
\hline
$\star$ & $3$&  \multicolumn{1}{|c|}{$5$ }\\
\cline{1-3}
$4$ & $6$\\
\cline{1-2}
$5$ & $7$\\
\cline{1-2}
\end{tabular}}\longrightarrow {\begin{tabular}{|c|c|c|}
\hline
 & $2$ & \multicolumn{1}{|c|}{$4$}\\
\hline
$3$ & $\star$ & \multicolumn{1}{|c|}{$5$ }\\
\cline{1-3}
$4$ & $6$\\
\cline{1-2}
$5$ & $7$\\
\cline{1-2}
\end{tabular}} \longrightarrow {\begin{tabular}{|c|c|c|}
\hline
 & $2$ & \multicolumn{1}{|c|}{$4$}\\
\hline
$3$ & \multicolumn{1}{|c|}{$5$} & $\star$\\
\cline{1-3}
$4$ & $6$\\
\cline{1-2}
$5$ & $7$\\
\cline{1-2}
\end{tabular}}\longrightarrow {\begin{tabular}{|c|c|c|}
\hline
 & $2$ & \multicolumn{1}{|c|}{$4$}\\
\hline
$3$ & \multicolumn{1}{|c|}{$5$ }\\
\cline{1-2}
$4$ & $6$\\
\cline{1-2}
$5$ & $7$\\
\cline{1-2}
\end{tabular}~=~T"}\\
T~&=~{\begin{tabular}{|c|c|c|}
\hline
 & $2$ & \multicolumn{1}{|c|}{$4$}\\
\hline
$\star$ & $3$ & \multicolumn{1}{|c|}{$6$ }\\
\cline{1-3}
$4$ & $5$\\
\cline{1-2}
$5$ & $7$\\
\cline{1-2}
\end{tabular}} \longrightarrow {\begin{tabular}{|c|c|c|}
\hline
 & $2$ & \multicolumn{1}{|c|}{$4$}\\
\hline
$3$ & $\star$ & \multicolumn{1}{|c|}{$6$ }\\
\cline{1-3}
$4$ & $5$\\
\cline{1-2}
$5$ & $7$\\
\cline{1-2}
\end{tabular}}
\longrightarrow {\begin{tabular}{|c|c|c|}
\hline
 & $2$ & \multicolumn{1}{|c|}{$4$}\\
\hline
$3$ & $5$ & \multicolumn{1}{|c|}{$6$}\\
\cline{1-3}
$4$ & $\star$\\
\cline{1-2}
$5$ & $7$\\
\cline{1-2}
\end{tabular}} \longrightarrow\\
&\longrightarrow {\begin{tabular}{|c|c|c|}
\hline
 & $2$ & \multicolumn{1}{|c|}{$4$}\\
\hline
$3$ & $5$ & \multicolumn{1}{|c|}{$6$}\\
\cline{1-3}
$4$ & $7$\\
\cline{1-2}
$5$ & $\star$\\
\cline{1-2}
\end{tabular}}
\longrightarrow {\begin{tabular}{|c|c|c|}
\hline & $2$ & \multicolumn{1}{|c|}{$4$}\\
\hline
$3$ & $5$ & \multicolumn{1}{|c|}{$6$}\\
\cline{1-3}
$4$ & $7$\\
\cline{1-2}
$5$\\
\cline{1-1}
\end{tabular}}~=~T".\endaligned
$$
\end{exe}

Appelons $S"$ le tableau de Young vide obtenu en supprimant la case point\'ee de $S$ et $\mu"=form(S")$. Le tableau $T"\setminus S"$ est encore semi standard. Si $(i,j)$ est le coin int\'erieur point\'e de $S$ et $(i",j")$ le coin ext\'erieur point\'e de $T"$, on pose $(T"\setminus S",(i",j"))=jdt(T\setminus S,(i,j))$. On peut inverser cette application.

Appelons inversion l'op\'eration qui consiste \`a prendre un tableau de Young semi standard $T\setminus S$ de forme $form(T\setminus S)=\lambda\setminus\mu$, \`a le plonger dans le plus petit rectangle le contenant (c'est \`a dire le rectangle de largeur $r$ et de hauteur $\ell_1$), puis \`a retourner ce rectangle et \`a remplacer chacune des entr\'ees $t_{ij}$ du tableau tordu ainsi obtenu par $n+1-t_{ij}$ et $\star$ par $\star$. Le tableau obtenu $T'\setminus S'=\sigma(T\setminus S)$ est encore un tableau semi standard tordu. Si on pointe un coin ext\'erieur de $T$, la case $\boxed{\star}$ est dans un coin int\'erieur de $S'$, et r\'eciproquement. Alors
$$
jdt^{-1}(T"\setminus S",(i",j"))=\sigma\circ jdt\circ \sigma(T"\setminus S",(i",j")).
$$
Par exemple le jeu de taquin appliqu\'e ci dessus s'inverse ainsi si $n=7$~:
$$\aligned
(T,(4,2))&={\begin{tabular}{|c|c|c|}
\hline
 & $2$ & \multicolumn{1}{|c|}{$4$}\\
\hline
$3$ & $5$ & \multicolumn{1}{|c|}{$6$ }\\
\cline{1-3}
$4$ & $7$\\
\cline{1-2}
$5$ & $\star$\\
\cline{1-2}
\end{tabular}}&\hfill
\sigma(T,(4,2))&={\begin{tabular}{|c|c|c|}
\hline
 & $\star$ & \multicolumn{1}{|c|}{$3$}\\
\hline
 & $1$ & \multicolumn{1}{|c|}{$4$}\\
\cline{1-3}
$2$ & $3$ & \multicolumn{1}{|c|}{$5$}\\
\cline{1-3}
$4$ & $6$\\
\cline{1-2}
\end{tabular}}\\
jdt\circ \sigma(T,(4,2))&={\begin{tabular}{|c|c|c|}
\hline
 & $1$ & \multicolumn{1}{|c|}{$3$}\\
\hline
 & $3$ & \multicolumn{1}{|c|}{$4$}\\
\cline{1-3}
$2$ & $5$ & \multicolumn{1}{|c|}{$\star$}\\
\cline{1-3}
$4$ & $6$\\
\cline{1-2}
\end{tabular}}&\hfill
\sigma\circ jdt\circ \sigma(T,(4,2))&=\begin{tabular}{|c|c|c|}
\hline
 & $2$ & \multicolumn{1}{|c|}{$4$}\\
\hline
$\star$ & $3$ & \multicolumn{1}{|c|}{$6$}\\
\cline{1-3}
$4$ & $5$\\
\cline{1-2}
$5$ & $7$\\
\cline{1-2}
\end{tabular}
\endaligned
$$

Le jeu de taquin est donc une application bijective
$$\aligned
jdt~:~&\bigcup_{\lambda\setminus\mu}SS(\lambda\setminus\mu)\times\{\text{coins int\'erieurs de }\mu\}~\longrightarrow\\
&\longrightarrow~\bigcup_{\lambda"\setminus\mu"}SS(\lambda"\setminus\mu")\times\{\text{coins ext\'erieurs de }\lambda"\}.
\endaligned
$$

Consid\`erons maintenant un tableau $T=(t_{ij})$ non quasi
standard et $s$ le plus grand entier tel que $T$ n'est pas quasi standard en $s$. Le haut de sa premi\`ere colonne est trivial :
$t_{s1}=s$, pour tout $j$, on a $t_{(s+1)j}<t_{s(j+1)}$ et $T$
poss\`ede une colonne de hauteur $s$.

On ajoute \`a ce tableau une colonne triviale, de hauteur $n-1$, dont on vide le sous tableau $S$ form\'e des $s$ cases sup\'erieures.

On pointe $S$ en son unique coin et on applique le jeu de taquin. La
case point\'ee se d\'eplace toujours vers la droite et 'sort' au
bout de la derni\`ere colonne de hauteur $s$. La ligne $s$ a juste
\'et\'e d\'ecal\'ee d'une case vers la gauche. On obtient un tableau
de premi\`ere colonne vide sur les $s-1$ premi\`eres cases et
triviale sur les $n-s$ cases restantes. On supprime cette colonne.
Si $s>1$, le tableau $T"$ obtenu n'est pas quasi standard en $s-1$, et peut \^etre en $s$, mais il est `quasi standard en tout $t>s$'.
On peut donc recommencer ce proc\'ed\'e et obtenir finalement un tableau quasi standard $T'$. Il est ais\'e de v\'erifier que cette proc\'edure r\'ealise une bijection entre l'ensemble des tableaux semi standards de forme $\lambda$ et l'union des tableaux quasi standards de forme plus petite que $\lambda$.
$$
SS^\lambda \longleftrightarrow \sqcup_{\mu \leq\lambda}~QS^{\mu}.
$$

D'autre part, on ordonne les tableaux de Young en disant que $T<S$ si $form(T)\leq form(S)$ et $form(T)\neq form(S)$ ou si $form(T)=form(S)$, mais qu'en lisant $T$ et $S$ colonne par colonne, de droite \`a gauche et de bas en haut, le premier couple d'entr\'ees diff\'erentes v\'erifie $s_{i,j}<t_{i,j}$.

Soit toujours $T$ un tableau non quasi standard en $s$ et quasi standard en tout $t>s$, de forme $\lambda$. Si la colonne num\'ero $\ell$ de $T$ a une hauteur sup\'erieure ou \'egale \`a $s$, et si $\partial^\ell T$ est le tableau obtenu en permutant les $s$ premi\`eres cases des colonnes num\'ero 1 et $\ell$ de $T$, on v\'erifie qu'en appliquant la relation de Pl\"ucker succesivement sur les colonnes num\'eros $i,i+1$, $1\leq i<\ell$, on obtient une relation
$$
T=\partial^\ell T+\sum_j S_j,
$$
o\`u $S_j<T$ pour tout $j$ (bien s\^ur $\partial^\ell T>T$). Lorsque la hauteur de la colonne num\'ero $\ell$ est $s$, on obtient une colonne triviale qui dispara\^\i t dans le quotient $\mathbb S^\bullet_{red}$. On appelle $(\partial^\ell T)'$ le tableau dans lequel on supprime cette colonne. On a $(\partial^\ell T)'<T$. On montre donc par r\'{e}currence que $\sqcup_{\mu \leq \lambda}~ QS^\mu$ est un syst\`eme g\'en\'erateur de ${\mathbb S}^{\lambda}_{\mathfrak n^+}$ dans ${\mathbb S}^{\bullet}_{red}$.\\
Comme ce syst\`eme de g\'en\'erateurs a pour cardinal la dimension
de ce module, on a prouv\'e :
\begin{theo} $($\cite{ABW}$)$

L'ensemble $QS^\bullet$ des tableaux quasi standards forme une base
de ${\mathbb S}^{\bullet}_{red}$, qui d\'{e}crit la stratification
de ce $\mathfrak n^+$-module ind\'ecomposable.

La r\'{e}union $\sqcup_{\mu\leq\lambda}~ QS^\mu$ forme une base de ${\mathbb S}^{\lambda}_{\mathfrak n^+}$.
\end{theo}

\section{Tableaux de Young semi standards symplectiques}

\

Cette section est consacr\'ee \`a rappeler la d\'efinition des tableaux de Young semi standards symplectiques. Cette notion a \'et\'e d\'evelopp\'ee en 1979 par C. De Concini (voir \cite{DeC}). En 1994, une autre description combinatoire des bases cristallines symplectiques, en termes de tableaux semi standards symplectiques, a \'et\'e pr\'esent\'ee par M. Kashiwara et T. Nakashima (voir \cite{KN}). En r\'ealit\'e, ces deux constructions sont \'equivalentes, une bijection explicite a \'et\'e donn\'ee par J. T. Sheats (\cite{SH}). Dans la suite, nous allons adapter la version des tableaux
semi standards symplectiques de De Concini en se r\'ef\'erant au travail de J. T. Sheats.

Rappelons aussi que R.G. Donnelly a donn\'e une construction explicite de l'action des \'el\'ements de $\mathfrak{sp}(2n)$ sur les bases de De Concini et de Kashiwara-Nakashima pour les repr\'esentations fondamentales $\mathbb S^{\langle\omega_k\rangle}$ (\cite{D}).

\subsection{Modules fondamentaux et colonnes symplectiques}

Utilisant l'ordre $1<2<\dots<n<\bar{n}<\dots<\bar{1}$, on \'{e}quipe $\mathbb{C}^{2n}$ de la base $(e_1,\dots,e_n, e_{\overline{n}},\dots,e_{\overline{1}})$ et de la forme symplectique
$$
\Omega=\displaystyle\sum e^\star_i \wedge e^\star_{\overline{i}}.
$$
Le groupe de Lie symplectique $SP(2n)$ est le groupe des matrices complexes $2n\times2n$ laissant $\Omega$ invariante. Son alg\`{e}bre de Lie $\mathfrak{sp}(2n)$ est simple de type $C_n$, c'est l'espace des matrices :
$$
X=\left(
\begin{array}{cc}
A & B\\
C & D\\
\end{array}\right),\quad A,~B,~C,~D\in Mat(n,n),~~ D=-^s A, B=^s B,C=^s C
$$
o\`u $^s$ est la sym\'etrie par rapport \`a la deuxi\`eme diagonale.

Une sous alg\`{e}bre de Cartan $\mathfrak h$ de $\mathfrak{sp}(2n)$ est la sous alg\`{e}bre des matrices diagonales $H=diag(\kappa_1,\dots,\kappa_n,-\kappa_n,\dots,-\kappa_1)$.
On pose $\theta_j(H)=\kappa_j$ et on choisit le syst\`eme de racines simples suivant~:
$$
\Delta=\{\alpha_i=\theta_i-\theta_{i+1},~~i=1, 2,\dots,n-1,~\alpha_n=2\theta_n\}.
$$
Remarquons que notre choix de racines simples est tel que, pour $\mathfrak{sp}(2n)$, la sous alg\`ebre $\mathfrak n^+=\sum_{\alpha>0}\mathfrak g^\alpha$ est l'espace des matrices de $\mathfrak{sp}(2n)$ qui sont strictement triangulaires sup\'erieures. On note $N^+$ le sous groupe analytique de $SP(2n)$ correspondant.

L'ensemble $\Lambda$ des poids entiers dominants est isomorphe \`a $\mathbb N^n$, en effet, $\lambda$ est entier dominant si et seulement si $\lambda=\sum_{k=1}^n a_k\omega_k$, o\`u les $a_k$ sont des entiers positifs ou nuls, et $\omega_k=\theta_1+\dots+\theta_k$ sont les poids fondamentaux de $\mathfrak{sp}(2n)$. Etudions d'abord ces modules fondamentaux.\\

\begin{theo} {\rm (\cite{FH})}

\

Si $k\geq2$, consid\'erons la fonction de contraction naturelle $\varphi_k$ d\'{e}finie par :

$$
{\varphi_k}(v_1\wedge\dots\wedge v_k)=\sum_{i<j}\Omega(v_i,v_j)(-1)^{i+j-1} v_1\wedge \dots \wedge \widehat{v_i}\wedge \dots \wedge \widehat{v_j} \wedge \dots \wedge v_k.
$$

Le noyau de $\varphi_k$ est un sous $\mathfrak{sp}(2n)$ module de $\wedge^k\mathbb C^{2n}$ isomorphe au module fondamental $\mathbb S^{\langle\omega_k\rangle}$ de $\mathfrak{sp}(2n)$ de plus haut poids $\omega_k$.
\end{theo}

Les $\mathfrak{sp}(2n)$ modules irr\'{e}ductibles fondamentaux $\mathbb S^{\langle\omega_k\rangle}$ sont ainsi r\'eali\-s\'es dans des sous espaces des $\mathfrak{sl}(2n)$ modules fondamentaux $\mathbb S^{\omega_k}$, pour $k=1,\dots,n$.

Tout $\mathfrak{sp}(2n)$ module simple $\mathbb{S}^{\langle\lambda\rangle}$ est le sous module du produit tensoriel
$$
Sym^{a_1}(\mathbb S^{\langle\omega_1\rangle})\otimes Sym^{a_2}(\mathbb S^{\langle\omega_2\rangle})\otimes\dots
\otimes Sym^{a_n}(\mathbb S^{\langle\omega_n\rangle})
$$
engendr\'e par le vecteur de plus haut poids.

Comme pour $SL(2n)$, consid\'erons les fonctions 'colonnes' suivantes d\'efinies sur $SP(2n)$~:
$$
\delta^{(k)}_{i_1,\dots,i_k}(g)=\langle e^\star_{i_1} \wedge \dots\wedge e^\star_{i_k}, ge_1\wedge \dots\wedge ge_k\rangle \hskip0.5cm (k \leq n,~~g\in SP(2n)).
$$
Ces fonctions ne sont pas ind\'ependantes. Par exemple, si $A,D\subset \{1,\dots,n\}$, si $A=\{p_1<p_2<\dots<p_s\}$, $D=\{q_1<\dots<q_t\}$, on pose
$$
e^{(\star)}_{A\overline{D}}=e^{(\star)}_{p_1}\wedge\dots\wedge e^{(\star)}_{p_s}\wedge e^{(\star)}_{\overline{q_t}}\wedge\dots\wedge e^{(\star)}_{\overline{q_1}}.
$$
Si $k=t+s+2\leq n$, on a
$$
\begin{aligned}
\langle e^\star_{A\overline{D}}\wedge \Omega, ge_{\{1,\dots,k\}}\rangle&=\sum_{i=1}^n\pm\langle e^\star_{A\cup\{i\}\overline{D\cup\{i\}}}, ge_{\{1,\dots,k\}}\rangle\\
&=\langle ^tg e^\star_{A\overline{D}}\wedge\Omega,e_{\{1\dots k\}}\rangle=0.
\end{aligned}
$$

Du th\'eor\`eme pr\'ec\'edent, on d\'eduit que ce sont les seules
relations homog\`enes de degr\'e 1 entre ces fonctions. On appellera
ces relations les relations de Pl\"ucker internes de
$\mathfrak{sp}(2n)$.

\begin{defi}

\

Soit $A, D \subset\{1,\dots,n\}$ tels que $k=\sharp A+ \sharp D \leq n$. Posons
$$
\begin{array}{c}
A\\
\overline{D}
\end{array}~=~\renewcommand{\arraystretch}{0.7}\begin{array}{|c|}
\hline
p_1\\
\hline
\vdots\\
\hline
p_s\\
\hline
\overline{q_t}\\
\hline
\vdots\\
\hline
\overline{q_1}\\
\hline
\end{array}~=~\delta^{(k)}_{p_1,\dots,p_s,\overline{q_t},\dots,\overline{q_1}}
$$
si $A=\{p_1<p_2<\dots<p_s\}$, $D=\{q_1<\dots<q_t\}$. Posons $I=A \cap D=\{i_1,\dots,i_r\}$.\\

On dit que la colonne est une colonne semi standard symplectique si $\{1,\dots,n\}\setminus A \cup D$ contient au moins un \'el\'ement $j>i_r$, deux \'el\'ements $j,j' >i_{r-1},$ etc\dots
\end{defi}

On montre $($\cite{DeC}$)$ que les colonnes semi standards symplectiques forment une base du module fondamental $\mathbb S^{\langle \omega_k\rangle}$.

On consid\`ere maintenant une colonne semi standard symplectique  $\renewcommand{\arraystretch}{0.7}\begin{array}{c}
A\\
\overline{D}\\
\end{array}$. On note $I=A \cap D$, $J=\{j_1,\dots,j_r\}$ la plus petite partie, pour l'ordre lexicographique, de $\{1,\dots,n\} \setminus A\cup D$ telle que $\# J= \# I$, $i_1<j_1$,\dots, $i_r<j_r$. On pose ($dble$ se lit `double') :
$$
C=(D\backslash I)\cup J,~~~~ B=(A\backslash I)\cup J,~~~~ dble\big(\renewcommand{\arraystretch}{0.7}\begin{array}{c}
A\\
\overline{D}
\end{array}\big)= \renewcommand{\arraystretch}{0.7}\begin{array}{cc}
A & B\\
\overline{C} & \overline{D}
\end{array}.
$$
Alors $dble\big(\renewcommand{\arraystretch}{0.7}\begin{array}{c}
A\\
\overline{D}
\end{array}\big)$ est un tableau de Young semi standard pour l'ordre choisi sur les indices : $1<2<\dots<n<\overline{n}<\dots<\overline{1}$.

\begin{exe} Supposons $n=4$, pour $\mathfrak{sp}(8)$, une colonne semi standard symplectique et son double est
$$
\begin{array}{c}
A\\
\overline{D}
\end{array}=\renewcommand{\arraystretch}{0.7}{\begin{array}{l}
\framebox{$1$}\\
\framebox{$2$}\\
\framebox{$\overline{1}$}\\
\end{array}}, ~~~~~dble\big(\renewcommand{\arraystretch}{0.7}\begin{array}{c}
A\\
\overline{D}
\end{array}\big)= \renewcommand{\arraystretch}{0.7}\begin{array}{cc}
A&B\\
\overline{C}&\overline{D}
\end{array}=\renewcommand{\arraystretch}{0.7}{\begin{array}{l}
\framebox{$1$}\framebox{$2$}\\
\framebox{$2$}\framebox{$3$}\\
\framebox{$\overline{3}$}\framebox{$\overline{1}$}\\
\end{array}}.
$$
\end{exe}
\


\subsection{Modules simples et tableaux semi standards symplectiques}

Soit $\lambda=\sum a_k\omega_k$ un poids entier dominant. Le module simple corres\-pondant $\mathbb S^{\langle\lambda\rangle}$ est le sous module engendr\'e par `le' vecteur de poids $\lambda$ dans
$$
Sym^{a_1}(\mathbb S^{\langle\omega_1\rangle})\otimes\dots\otimes Sym^{a_n}(\mathbb S^{\langle\omega_n\rangle}).
$$

Il est donc engendr\'e par les tableaux de Young de forme
$$
\lambda= (a_1,\dots,a_n,0,\dots,0)
$$
dont toutes les colonnes sont semi standards symplectiques. Une base
de $\mathbb S^{\langle\lambda\rangle}$ a \'et\'e d\'etermin\'ee par
G. de Concini (\cite{DeC}).

\begin{defi}

\

Soit $T$ un tableau de forme $\lambda$ dont toutes les colonnes sont semi standards symplectiques. Le tableau $dble(T)$ est le tableau obtenu en juxtaposant les doubles des colonnes de $T$.

On dit que $T$ est un tableau semi standard symplectique (ou semi
standard pour $\mathfrak{sp}(2n)$) si $dble(T)$ est un tableau semi
standard (pour $\mathfrak{sl}(2n)$).
\end{defi}

Alors

\begin{theo} (\cite{DeC})

L'ensemble $SS^{\langle\lambda\rangle}$ des tableaux de Young semi
standards symplectiques de forme $\lambda$ est une base du
$\mathfrak{sp}(2n)$ module simple
$\mathbb{S}^{\langle\lambda\rangle}$.
\end{theo}

\begin{exe}

\

Pour $n=3$ (cas de $\mathfrak{sp}(6)$), et $\lambda=\omega_2+\omega_3$, le tableau suivant est semi standard symplectique~:
$$
\begin{tabular}{|c|c|}
\hline

\raisebox{-2pt}{$1$} & \raisebox{-2pt}{$2$}\\
\hline
\raisebox{-2pt}{$2$} & \raisebox{-2pt}{$\overline{2}$}\\
\cline{1-2}
\raisebox{-2pt}{$\overline{2}$}\\
\cline{1-1}
\end{tabular}~~ \text{ en effet } ~~dble~(T)={\begin{tabular}{|c|c|c|c|}
\hline
\raisebox{-2pt}{$1$} & \raisebox{-2pt}{$2$} & \raisebox{-2pt}{$2$} & \raisebox{-2pt}{$3$}\\
\hline
\raisebox{-2pt}{$2$} & \raisebox{-2pt}{$3$} & \raisebox{-2pt}{$\overline{3}$} & \raisebox{-2pt}{$\overline{2}$}\\
\cline{1-4}
\raisebox{-2pt}{$\overline{3}$} & \raisebox{-2pt}{$\overline{2}$}\\
\cline{1-2}
\end{tabular}}.
$$
\end{exe}


\section{Alg\`ebre de forme et alg\`ebre de forme r\'eduite}

\

Consid\'erons la somme de tous les $\mathfrak{sp}(2n)$ modules simples :
$$
\mathbb{S}^{\langle\bullet\rangle}=\bigoplus_{\lambda\in\Lambda}~\mathbb{S}^{\langle\lambda\rangle}.
$$

L'ensemble $SS^{\langle\bullet\rangle}$ de tous les tableaux semi standards symplectiques est donc une base de ce module.\\

Comme pour toute alg\`ebre de Lie semi simple, cette somme peut \^etre munie d'une multiplication qui en fait une alg\`ebre commutative. Pour $\mathfrak{sp}(2n)$, on peut r\'ealiser cette structure explicitement, exactement comme pour $\mathfrak{sl}(2n)$.\\

Notons $\mathbb C[SP(2n)]^{N^+}$ l'espace des fonctions polynomiales sur $SP(2n)$ qui sont invariantes par multiplications \`a droite par les matrices de $N^+$. C'est un $SP(2n)$ module pour l'action \`a gauche~:
$$
(g.f)(g_1)=f(^tgg_1)\qquad (g,~g_1\in SP(2n),~~f\in\mathbb C[SP(2n)]^{N^+}).
$$
Comme c'est aussi une somme de modules de dimension finie, il se d\'ecompose en somme de modules irr\'eductibles $\mathbb S^{\langle \lambda \rangle}$. Ses vecteurs de poids dominant $f^\lambda$ sont des fonctions polynomiales invariantes sous la multiplication \`a droite par $N^+$ et \`a gauche par $^tN^+$. Par la m\'ethode du pivot de Gauss, ces fonctions sont caract\'eris\'ees par leur valeur sur les matrices diagonales de $SP(2n)$ :
$$\aligned
f^\lambda(g)&=f^\lambda\left(diag(\delta^{(1)}_1(g),\dots,\frac{\delta^{(n)}_{1\dots n}(g)}{\delta^{(n-1)}_{1\dots(n-1)}(g)},\frac{\delta^{(n-1)}_{1\dots(n-1)}(g)}{\delta^{(n)}_{1\dots n}(g)},\dots,\frac{1}{\delta^{(1)}_1(g)})\right)\\
&=\sum_{p_j\in\mathbb Z} c_{p_1,\dots,p_n}\left(\delta^{(1)}_1(g)\right)^{p_1}\dots\left(\frac{\delta^{(n-1)}_{1\dots(n-1)}(g)}{\delta^{(n)}_{1\dots n}(g)}\right)^{p_n}
\endaligned
$$
(la derni\`ere somme est finie). En faisant agir $\mathfrak h$ sur cette fonction, on voit que la somme ne contient qu'un terme et que
$$
\lambda=(p_1-p_2)\omega_1+(p_2-p_3)\omega_2+\dots+(p_{n-1}-p_n)\omega_{n-1}+p_n\omega_n.
$$
Comme $\lambda$ est dominant entier, les $p_k$ sont entiers et v\'erifient $p_1\geq p_2\geq\dots\geq p_n\geq0$. Pour chaque $\lambda$ de $\Lambda$, l'espace des fonctions $f^\lambda$ est de dimension 1, ou :
$$
\mathbb C[SP(2n)]^{N^+}\simeq\bigoplus_{\lambda\in\Lambda}\mathbb S^{\langle \lambda \rangle}=\mathbb S^{\langle \bullet \rangle}.
$$

Cette identification fait de $\mathbb S^{\langle \bullet \rangle}$
une alg\`ebre commutative.

\begin{defi}

\

On appelle alg\`ebre de forme de $\mathfrak{sp}(2n)$ l'alg\`ebre
$\mathbb S^{\langle \bullet \rangle}$ munie de la multiplication
d\'efinie ci-dessus.
\end{defi}

Gr\^ace aux r\'esultats pr\'ec\'edents, on a:

\begin{prop}

\

L'alg\`ebre de forme $\mathbb S^{\langle \bullet \rangle}$ de
$\mathfrak{sp}(2n)$ est le quotient de la sous alg\`ebre $\mathbb
S^{\bullet}_{(n)} =\oplus_{\lambda=(a_1,\dots,a_n,0,\dots,0)}\mathbb
S^\lambda$ de l'alg\`ebre de forme de $\mathfrak{sl}(2n)$ par
l'id\'{e}al $J^{\langle \bullet \rangle}$ engendr\'e par les
relations de Pl\"ucker internes.
\end{prop}

On peut donc \'ecrire :
$$
\mathbb S^{\langle \bullet \rangle}\simeq \mathbb{C}{[SP(2n)]}^{N^+} \simeq \mathbb{C}{[\delta^{(r)}_{i_1,\dots,i_r},~~r\leq n]}/ \mathcal{PL}
$$
o\`u $\mathcal{PL}$ est l'id\'eal des relations de Pl\"ucker externes sur les couples de colonnes de hauteur $\leq n$ (relations homog\'enes de degr\'ee deux) et des relations de Pl\"ucker internes (homog\`enes de degr\'e un).\\

Comme pour $SL(n)$, en restreignant les fonctions
$\delta^{(r)}_{i_1,\dots,i_r}$ ($r\leq n$) \`a $N^-=^t N^+$, on
d\'{e}finit l'alg\`{e}bre forme r\'{e}duite pour
$\mathfrak{sp}(2n)$.

\begin{defi}

\

On appelle alg\`ebre de forme r\'eduite de $\mathfrak{sp}(2n)$ et on note $\mathbb{S}^{\langle \bullet \rangle}_{red}$ le quotient :
$$
\mathbb S^{\langle \bullet\rangle}_{red}=\mathbb S^{\langle \bullet \rangle}\big{/}<\delta^{(k)}_{1,\dots,k}-1>,\hskip0.3cm k=1,2,\dots,n.
$$
\end{defi}

\begin{theo}

\

\begin{itemize}
\item[i)] $\mathbb S^{\langle \bullet \rangle}_{red}$ est un $\mathfrak n^+$ module ind\'ecomposable.\\
\item[ii)] $\mathbb S^{\langle \bullet\rangle}_{red}$ est l'union des $\mathbb S^{\langle\lambda\rangle}_{\mathfrak n^+}$, stratifi\'ee par :
$$
\mu \leq \lambda ~~\Longleftrightarrow ~~\mathbb S^{\langle \mu\rangle}_{\mathfrak n^+} \subset \mathbb S^{\langle \lambda \rangle}_{\mathfrak n^+}.
$$

\item[iii)] Tout $\mathfrak n^+$ module monog\`ene localement nilpotent est un quotient d'un des $\mathbb{S}^{\langle \lambda \rangle}_{\mathfrak n^+}$.\\
\item[iv)] On a $\mathbb S^{\langle \bullet \rangle}_{red}=\mathbb S^\bullet_{(n)~~red} \big/ J^{\langle \bullet \rangle}$ o\`u
$$
\mathbb S^\bullet_{(n)~~red}=\bigoplus_{\lambda=(\lambda_1,\dots,\lambda_n,0,\dots,0)} \mathbb S^\lambda\big{/}< \delta^{(k)}_{1,\dots,k}-1,~~k\leq n>.
$$
\end{itemize}
\end{theo}
{\bf{Preuve:}}

\

Les preuves de i), ii) et iii) sont identiques \`a celles de \cite{ABW} pour le cas de $\mathfrak{sl}(n)$.\\

\noindent
iv) Le diagramme suivant
$$
\begin{array}{ccccc}
 & & \phi & & \\
\mathbb S^\bullet_{(n)} & & \longrightarrow & \mathbb S^{\langle \bullet\rangle}=\mathbb S^\bullet_{(n)} {_{\big /J^{\langle\bullet\rangle}}} & \\
 & & & & \\
\pi_1 \downarrow & & & \downarrow \pi & \\
 & & \phi_1 & & \\
\mathbb S_{(n)~~red}^\bullet=\mathbb S_{(n)}^\bullet\;_{ \big/ < \delta^{(k)}_{1,\dots,k}-1,~~k\leq n>} & & \longrightarrow & \mathbb S_{red}^{\langle \bullet \rangle}=\mathbb S^{\langle \bullet\rangle}~_{ \big/< \delta^{(k)}_{1,\dots,k}-1,~~k\leq n>} & \\
 & & & \hskip1.2cm \simeq \mathbb S_{red}^\bullet {_{\big /J^{\langle\bullet\rangle}}} &
\end{array}
$$
est commutatif c'est \`a dire $\pi \circ  \phi= \phi_1 \circ\pi_1$.
On a donc bien $\mathbb S_{red}^{\langle \bullet\rangle}= \mathbb
S_{(n)~~red}^\bullet{_{\big/ J^{\langle\bullet\rangle}}}$.

$\hfill\square$

\section{Tableaux de Young quasi standards symplectiques}


\

A partir de maintenant, on notera aussi $f(A,D)$ la colonne semi standard symplectique $\begin{array}{c}
A\\
\overline{D}\\
\end{array}$. Rappellons nos notations $I=A\cap D$, $J$ est la plus petite partie `\`a droite de $I$' dans le compl\'ementaire de $A\cup D$ et $B=(A\setminus I)\cup J$, $C=(D\setminus I)\cup J$.\\

Inversement, si $B$ et $C$ sont connus, on peut retrouver $I$, $J$, $A$ et $D$. En effet on a alors $J=B\cap C$ et $I$ est la plus grande partie \`a gauche de $J$ dans le compl\'ementaire de $B\cup C$, ayant le m\^eme nombre d'\'el\'ements que $J$.\\

Soit donc $B$ et $C$ deux parties de $\{1,\dots,n\}$. Posons $J=B\cap C=\{j_1<\dots<j_r\}$ et d\'efinissons $I$ comme la plus grande partie $\{i_1<i_2<\dots<i_r\}$, pour l'ordre lexicographique, de $\mathbb Z\setminus(B\cup C)$ telle que $i_k<j_k$ pour tout $k$. On pose enfin :
$$
A=(B\setminus J)\cup I,\quad D=(C\setminus J)\cup I\quad\hbox{ et }\quad g(B,C)=f(A,D)=\begin{array}{c}
A\\
\overline D\\
\end{array}.
$$

\begin{defi}

\

Soit $T$ un tableau de Young semi standard symplectique. Nous dirons
que $T$ est quasi standard symplectique si $dble(T)$ est quasi
standard (pour $\mathfrak{sl}(2n)$).
\end{defi}

Notons $SS^{\langle \lambda \rangle}$ l'ensemble des tabeaux semi standards symplectiques de forme $\lambda=(a_1,a_2,\dots,a_n)$ ayant $a_1$ colonnes de hauteur 1, \dots, $a_n$ colonnes de hauteur $n$. De m\^eme, notons $QS^{\langle \lambda\rangle}$ l'ensemble des tableaux quasi standards symplectiques de forme $\lambda$ et $NQS^{\langle\lambda\rangle}$ l'ensemble des tableaux semi standards non quasi standards symplectiques de forme $\lambda$.\\

Remarquons qu'un tableau de Young $T$ peut \^etre quasi standard pour $\mathfrak{sl}(2n)$ sans que son double le soit. En voici un exemple
$$
T=\begin{tabular}{|c|c|}
\hline
\raisebox{-2pt}{$1$} & \raisebox{-2pt}{$2$}\\
\hline
\raisebox{-2pt}{$2$} & \raisebox{-2pt}{$\overline 2$}\\
\cline{1-2}
\raisebox{-2pt}{$\overline 2$}\\
\cline{1-1}
\end{tabular}~~ \Longrightarrow ~~dble(T)=\begin{tabular}{|c|c|c|c|}
\hline
\raisebox{-2pt}{$1$} & \raisebox{-2pt}{$1$} & \raisebox{-2pt}{$2$} & \raisebox{-2pt}{$3$}\\
\hline
\raisebox{-2pt}{$2$} & \raisebox{-2pt}{$3$} & \raisebox{-2pt}{$\overline{3}$} & \raisebox{-2pt}{$\overline{2}$}\\
\cline{1-4}
\raisebox{-2pt}{$\overline{3}$} & \raisebox{-2pt}{$\overline{2}$}\\
\cline{1-2}
\end{tabular}
$$
$T$ est quasi standard mais $dble(T)$ ne l'est pas.\\

On dira qu'un tableau $T$ semi standard symplectique est poussable en $s$, et on notera $T\in NQS_s$ si $dble(T)=(t_{i,j})$ a la propri\'et\'e $t_{s,j+1}<t_{s+1,j}$, pour  tout $j$ pour lesquels ces deux entr\'ees existent. On remarque d'abord que chaque colonne de $T$, \'el\'ement de $NQS_s$ se d\'ecompose.

\begin{lem} Soit $T$ un tableau de $NQS_s$, $c=\begin{array}{c}A\\ \overline{D}\end{array}$ une colonne de $T$ et $dble(c)=\begin{array}{cc}A&B\\ \overline{C}&\overline{D}\end{array}$ son double. Supposons $s\leq \#A$. Soit $\alpha$ un nombre entier tel que $b_s\leq\alpha<a_{s+1}$ (si $s=\#A$, on choisit $\alpha\geq b_s$ seulement). Pour toute partie $X$ de $[1,n]$, on pose $X^{\leq\alpha}=X\cap[1,\alpha]$ et $X^{>\alpha}=X\cap]\alpha,n]$. Alors
\begin{itemize}
\item[1.] La colonne $\begin{array}{c}A^{\leq\alpha}\\ \overline{D}^{\leq\alpha}\end{array}$ est semi standard pour $\mathfrak{sp}(2\alpha)$ et son double est $\begin{array}{cc}A^{\leq\alpha}&B^{\leq \alpha}\\ \overline{C}^{\leq\alpha}&\overline{D}^{\leq\alpha}\end{array}$.\\
\item[2.] La colonne $\begin{array}{c}A^{>\alpha}\\ \overline{D}^{>\alpha}\end{array}$, index\'ee par $[\alpha+1,n]\cup[\overline{n},\overline{\alpha+1}]$, est semi standard pour $\mathfrak{sp}(2(n-\alpha))=\mathfrak{sp}(2]\alpha,n])$ et son double est $\begin{array}{cc}A^{>\alpha}&B^{>\alpha}\\ \overline{C}^{>\alpha}&\overline{D}^{>\alpha}\end{array}$.\\
\end{itemize}
\end{lem}

\noindent
{\bf Preuve}

Par hypoth\`ese, $b_s\leq\alpha<a_{s+1}$ $(b_s\leq\alpha$ si
$s=\#A)$ et $I^{\leq\alpha}=A^{\leq\alpha}\cap D^{\leq\alpha}$, donc
les \'el\'ements $j_1,\dots,j_r$ qui sont dans le compl\'ementaire
de $(A\cup D)$ et dans $[1,\alpha]$ comprennent les \'el\'ements de
$J\cap([1,\alpha]\setminus (A^{\leq\alpha}\cup D^{\leq\alpha})$, donc il y en a suffisamment pour que la colonne
$\begin{array}{c}A^{\leq \alpha}\\
\overline{D}^{\leq\alpha}\end{array}$ soit semi standard pour
$\mathfrak{sp}(2\alpha)$. Par construction, $J^{\leq\alpha}$ est la
plus petite partie de $[1,\alpha]\setminus(A^{\leq\alpha}\cup
D^{\leq\alpha})$, de cardinal $\#I^{\leq\alpha}$ et qui contienne 1
\'el\'ement plus grand que le premier \'el\'ement de
$I^{\leq\alpha}$, un deuxi\`eme \'el\'ement plus grand que le
deuxi\`eme \'el\'ement de $I^{\leq\alpha}$, etc \dots, donc
$$
dble(\begin{array}{c}A^{\leq\alpha}\\ \overline{D}^{\leq\alpha}\end{array})=\begin{array}{cc}A^{\leq\alpha}&B^{\leq\alpha}\\ \overline{C}^{\leq\alpha}&\overline{D}^{\leq\alpha}\end{array}.
$$

On en d\'eduit que si $I^{>\alpha}=I\setminus I^{\leq\alpha}$ et
$J^{>\alpha}=J\setminus J^{\leq\alpha}$, alors $J^{>\alpha}$ a le
m\^eme cardinal que $I^{>\alpha}$, chaque \'el\'ement de
$I^{>\alpha}$ est major\'e par un \'el\'ement de $J^{>\alpha}$ et
$J^{>\alpha}$ est la plus petite partie de $[\alpha+1,n]\setminus
A^{>\alpha}\cup B^{>\alpha}$ ayant cette propri\'et\'e. Ceci
ach\`eve la preuve du lemme.

$\hfill\square$

On notera une telle colonne
$$
c=f^{\leq\alpha}(A^{\leq\alpha},D^{\leq\alpha})\uplus f^{>\alpha}(A^{>\alpha},D^{>\alpha})=\begin{array}{ccc}
&A^{\leq\alpha}&\\
\hline
\hline
&A^{>\alpha}&\\
&\overline{D}^{>\alpha}&\\
&\overline{D}^{\leq\alpha}&
\end{array}.
$$
On notera aussi :
$$\aligned
c&=g^{\leq\alpha}(B^{\leq\alpha},C^{\leq\alpha})\uplus g^{>\alpha}(B^{>\alpha},C^{>\alpha}),\\
dble(c)&=\begin{array}{cccc}
&A^{\leq\alpha}&B^{\leq\alpha}&\\
\hline
\hline
&A^{>\alpha}&B^{>\alpha}&\\
&\overline{C}^{>\alpha}&\overline{D}^{>\alpha}&\\
&\overline{C}^{\leq\alpha}&\overline{D}^{\leq\alpha}&
\end{array}.
\endaligned
$$

De m\^eme, si $c$ est une colonne de hauteur $k$ qui n'est pas quasi standard en $s$ et $s>\#A$, alors $k-s+1=\#A+\#D-s+1\leq\#D$. On choisit $\alpha$ tel que $c_{k-s}<\alpha\leq d_{k-s+1}$, la colonne $f^{<\alpha}(A^{<\alpha},D^{<\alpha})$ est semi standard pour $\mathfrak{sp}(2(\alpha-1))$, et la colonne $f^{\geq\alpha}(A^{\geq\alpha},D^{\geq\alpha})$ est semi standard pour $\mathfrak{sp}(2(n-\alpha+1))$, index\'ee par $[\alpha,n]$. Le double de ces colonnes est comme ci-dessus. On notera la colonne $c$:
$$
c=f^{<\alpha}(A^{<\alpha},D^{<\alpha})\uplus f^{\geq\alpha}(A^{\geq\alpha},D^{\geq\alpha})=\begin{array}{ccc}
&A^{<\alpha}&\\
&A^{\geq\alpha}&\\
&\overline{D}^{\geq\alpha}&\\
\hline
\hline
&\overline{D}^{<\alpha}&
\end{array}.
$$
On notera aussi :
$$\aligned
c&=g^{<\alpha}(B^{<\alpha},C^{<\alpha})\uplus g^{\geq\alpha}(B^{\geq\alpha},C^{\geq\alpha}),\\
dble(c)&=\begin{array}{cccc}
&A^{<\alpha}&B^{<\alpha}&\\
&A^{\geq\alpha}&B^{\geq\alpha}&\\
&\overline{C}^{\geq\alpha}&\overline{D}^{\geq\alpha}&\\
\hline
\hline
&\overline{C}^{<\alpha}&\overline{D}^{<\alpha}&
\end{array}.
\endaligned
$$

\begin{lem}
Avec les hypoth\`eses du lemme pr\'ec\'edent, on pose $A'=A\setminus\{a_s\}$ et $f(A',D)=g(B',C')$ (dans $\mathfrak{sp}(2n)$). Alors
$$
b_s\notin A'\cup B'\cup C'\cup D.
$$
Ou
$$\aligned
f(A',D)&=f^{<b_s}(A'^{<b_s},D^{<b_s})\uplus f^{>b_s}(A^{>b_s},D^{>b_s})\\
&=g^{<b_s}(B'^{<b_s},C'^{<b_s})\uplus g^{>b_s}(B^{>b_s},C^{>b_s}).
\endaligned
$$
\end{lem}

\noindent
{\bf Preuve}

\begin{itemize}
\item[{\bf Cas 1}] $a_s=b_s$\\ Comme $f^{\leq b_s}(A^{\leq b_s},D^{\leq b_s})$ est une colonne symplectique, alors $a_s$ n'appartient pas \`a $I^{\leq b_s}$. Donc $I^{<b_s}=I^{\leq b_s}=I'^{\leq b_s}$ et $J^{<b_s}=J^{\leq b_s}=J'^{\leq b_s}$. Par suite $b_s$ n'appartient pas \`a $B'^{\leq b_s}$, $b_s\notin B'$, on a $b_s\notin A'\cup D$ et $b_s\notin C'$.
\item[{\bf Cas 2}] $a_s<b_s$\\ Cela veut dire $b_s\in J^{\leq b_s}$, plus pr\'ecis\'ement
$$
I^{\leq b_s}=\{a_{t_1}<\dots<a_{t_r}\}\quad\text{ et }\quad J^{\leq b_s}=\{b_{u_1}<\dots<b_{u_r}=b_s\}.
$$
On a deux sous-cas :
$\bullet$ Si $t_r=s$, en enlevant $a_{t_r}=a_s$ de $A$, on enl\`eve $a_s$ de $I^{\leq b_s}$~: $I'^{\leq b_s}=I^{\leq b_s}\setminus\{a_s\}$ et donc $b_s$ de $J'^{\leq b_s}$, c'est \`a dire $b_s\notin A'\cup B'\cup C'\cup D$.\\
$\bullet$ Si $t_r<s$, alors par construction, $a_s\notin I^{\leq b_s}$, en supprimant $a_s$ pour construire $A'$, on a $a_s\in[1,b_s]\setminus (A'\cup D)$ et donc $J^{\leq b_s}$ devient $J'^{\leq b_s}=\{b_{u_1},\dots,b_{u_{r-1}},a_s\}$ (l'ordre n'est peut \^etre pas pr\'eserv\'e). Donc $b_s\notin A'\cup B'\cup C'\cup D$.\\
\end{itemize}
Finalement puisque $f^{>b_s}(A^{>b_s},D^{>b_s})=g^{>b_s}(B^{>b_s},C^{>b_s})$ et que $b_s$ n'appartient pas \`a $D$, on a $A^{<a_s}=A'^{<b_s}$ et la derni\`ere relation :
$$\aligned
f(A',D)&=f^{<b_s}(A'^{<b_s},D^{<b_s})\uplus f^{>b_s}(A^{>b_s},D^{>b_s})\\
&=g^{<b_s}(B'^{<b_s},C'^{<b_s})\uplus g^{>b_s}(B^{>b_s},C^{>b_s})\\
&=g(B',C').
\endaligned
$$
$\hfill\square$

Si $s>\#A$, posons $t=k-s+1=\#A+\#D-s+1$, alors
$$
f(A,D)=g(B,C)=g^{<d_t}(B^{<d_t},C^{<d_t})\uplus g^{\geq d_t}(B^{\geq d_t},C^{\geq d_t}),
$$
si $C'=C\setminus\{c_t\}$, et $g(B,C')=f(A',D')$, alors $d_t$ n'appartient pas \`a $B\cup C'\cup A'\cup D'$, $C^{>d_t}=C'^{>d_t}$ et
$$\aligned
g(B,C')&=g^{<d_t}(B^{<d_t},C^{<d_t})\uplus g^{>d_t}(B^{>d_t},C'^{>d_t})\\
&=f^{<d_t}(A^{<d_t},D^{<d_t})\uplus f^{>d_t}(A'^{>d_t},D'^{>d_t}).
\endaligned
$$

\begin{lem}
Soit $f(A',D)$ la colonne du lemme pr\'ec\'edent. On garde les notations de ce lemme.
\begin{itemize}
\item Soit $u\geq b_s$, $u\notin B^{>b_s}$, $B"=B'\cup\{u\}$ et $C"=C'$, alors:
$$\aligned
g(B",C")&=g^{<b_s}(B'^{<b_s},C'^{<b_s})\uplus g^{\geq b_s}(\{u\}\cup B^{>b_s},C^{>b_s})\\
&=f^{<b_s}(A^{<b_s},D^{<b_s})\uplus f^{\geq b_s}(A"^{\geq b_s},D"^{\geq b_s})\\
&=f(A",D").
\endaligned
$$
\item Supposons que $A^{>b_s}=\emptyset$. Soit $v\geq b_s$, $v\notin D^{>b_s}$ soit $A"=A'$ et $D"=D'\cup\{v\}$, alors~:
$$\aligned
f(A",D")&=f^{<b_s}(A^{<b_s},D^{<b_s})\uplus f^{\geq b_s}(\emptyset,D^{>b_s}\cup\{v\})\\
&=g^{<b_s}(B'^{<b_s},C'^{<b_s})\uplus g^{\geq b_s}(\emptyset,D^{>b_s}\cup\{v\})\\
&=g(B",C").
\endaligned
$$
\item Supposons que $A^{>b_s}=D^{>b_s}=\emptyset$. Soit $v< b_s$, $v\notin D^{<b_s}$, soit $A"=A'$ et $D"=D\cup\{v\}$, alors~:
$$\aligned
f(A",D")&=f^{\leq b_s}(A^{<b_s},D^{<b_s}\cup\{v\})\\
&=g^{\leq b_s}(B"^{\leq b_s},C"^{\leq b_s})\\
&=g(B",C").
\endaligned
$$
\end{itemize}

\end{lem}

\noindent
{\bf Preuve}

Dans le premier cas, puisque la colonne $g^{>b_s}(B^{>b_s},C^{>b_s})$ est semi standard pour $\mathfrak{sp}(2]b_s,n])$, la colonne $g^{\geq b_s}(\{u\}\cup B^{>b_s},C^{>b_s})$ est semi standard pour $\mathfrak{sp}(2[b_s,n])$. En effet, si $u>b_s$ l'ensemble $J'^{>b_s}$ devient $J"^{> b_s}=J'^{>b_s}$ ou $J"^{>b_s}=J'^{>b_s}\cup\{u\}$, mais dans ce dernier cas, l'ajout de l'indice $b_s$ garantit que la colonne reste semi standard. Si $u=b_s$, alors $J"^{\geq b_s}=J^{>b_s}$ et on n'a pas besoin de l'indice $b_s$ pour construire $I"^{\geq b_s}=I'^{>b_s}$. Ceci prouve le premier cas.\\

Dans le cas 2, si $v>b_s$, il n'y a rien \`a prouver, si $v=b_s$, on doit prendre $f^{\geq b_s}$ \`a la place de $f^{>b_s}$.\\

Dans le cas 3, $f(A",D")$ est semi standard puisqu'on ajoute, en
m\^eme temps que $v$, un indice ($b_s$) apr\`es les indices de la
colonne semi standard $f^{< b_s}(A^{<b_s},D^{<b_s})$ de
$\mathfrak{sp}(2[1,b_s[)$.

$\hfill\square$

Si maintenant $s>\#A$, on pose comme plus haut $t=k-s+1$ et on
consid\`ere la colonne
$$
g(B,C')=g^{<d_t}(B^{<d_t},C^{<d_t})\uplus g^{>d_t}(B^{>d_t},C^{d_t})=f(A',D').
$$
Soit $v\leq d_t$, $v\notin D^{<d_t}$, et $A"=A'$, $D"=D'\cup\{v\}$ alors le m\^eme argument que le cas 1 ci dessus donne~:
$$\aligned
f(A",D")&=f^{<d_t}(A^{<d_t},D^{<d_t})\uplus f^{\geq d_t}(A'^{>d_t},D'^{>d_t}\cup\{v\})\\
&=f^{\leq d_t}(A"^{\leq d_t},D"^{\leq d_t})\uplus f^{>d_t}(A'^{>d_t},D'^{>d_t})\\
&=g^{\leq d_t}(B"^{\leq d_t},C"^{\leq d_t})\uplus g^{>d_t}(B'^{>d_t},C'^{>d_t})\\
&=g(B",C").
\endaligned
$$

\begin{exe}

\

On consid\'{e}re le tableau dans $\mathfrak{sp}(18)$:
$$T=\begin{tabular}{|c|c|c|c|} \hline
\raisebox{-2pt}{$1$}& \raisebox{-2pt}{$4$} & \raisebox{-2pt}{$5$} & \raisebox{-2pt}{$7$}\\
\hline
  \raisebox{-2pt}{$2$}& \raisebox{-2pt}{$5$} & \raisebox{-2pt}{$7$}&  \raisebox{-2pt}{$\overline{7}$}\\
\hline
 \raisebox{-2pt}{$3$} & \raisebox{-2pt}{$6$} & \raisebox{-2pt}{$\overline{5}$}& \raisebox{-2pt}{$\overline{4}$}\\
\hline
\raisebox{-2pt}{$7$} & \raisebox{-2pt}{$\overline{5}$}& \raisebox{-2pt}{$\overline{3}$}\\
\cline{1-3}
\raisebox{-2pt}{$8$}\\
\cline{1-1}
\raisebox{-2pt}{$\overline{8}$}\\
\cline{1-1}
\raisebox{-2pt}{$\overline{3}$}\\
\cline{1-1}
\raisebox{-2pt}{$\overline{2}$}\\
\cline{1-1}\raisebox{-2pt}{$\overline{1}$}\\
\cline{1-1}
\end{tabular} \in NQS_3.$$
On note $c= \begin{array}{c}
A\\
\overline{D}\\
\end{array}$ la premi\`{e}re colonne de $T$, alors $$\aligned dble(c)=& \begin{array}{cccc}
&A&B&\\
&\overline{C}&\overline{D}&
\end{array}
=\begin{tabular}{|c|c|} \hline
\raisebox{-2pt}{$1$}& \raisebox{-2pt}{$4$} \\
\hline
  \raisebox{-2pt}{$2$}& \raisebox{-2pt}{$5$}\\
\hline
 \raisebox{-2pt}{$3$} & \raisebox{-2pt}{$6$}\\
\hline
\raisebox{-2pt}{$7$} & \raisebox{-2pt}{$7$}\\
\hline
\raisebox{-2pt}{$8$} & \raisebox{-2pt}{$9$}\\
\hline
\raisebox{-2pt}{$\overline{9}$}& \raisebox{-2pt}{$\overline{8}$}\\
\cline{1-2}
\raisebox{-2pt}{$\overline{6}$}& \raisebox{-2pt}{$\overline{3}$}\\
\cline{1-2}
\raisebox{-2pt}{$\overline{5}$}&\raisebox{-2pt}{$\overline{2}$}\\
\cline{1-2}
\raisebox{-2pt}{$\overline{4}$}&\raisebox{-2pt}{$\overline{1}$}\\
\cline{1-2}
\end{tabular}\endaligned.$$
Avec les notations des lemmes pr\'{e}c\'{e}dents, on a :
$s=3<\#A$ et $ b_s = \alpha= 6 <a_{s+1}$.\\
La colonne $\begin{array}{c}A^{\leq\alpha}\\
\overline{D}^{\leq\alpha}\end{array}= \begin{tabular}{|c|} \hline
\raisebox{-2pt}{$1$}\\
\hline
  \raisebox{-2pt}{$2$}\\
\hline
 \raisebox{-2pt}{$3$} \\
\cline{1-1}
\raisebox{-2pt}{$\overline{3}$}\\
\cline{1-1}
\raisebox{-2pt}{$\overline{2}$}\\
\cline{1-1}\raisebox{-2pt}{$\overline{1}$}\\
\cline{1-1}
\end{tabular}$ est semi standard pour
$\mathfrak{sp}(12)$ et $\begin{array}{c}A^{>\alpha}\\
\overline{D}^{>\alpha}\end{array}=\begin{tabular}{|c|} \hline
\raisebox{-2pt}{$7$}\\ \hline
  \raisebox{-2pt}{$8$}\\
\cline{1-1}
\raisebox{-2pt}{$\overline{8}$}\\\cline{1-1}
\end{tabular}$ est semi standard pour $\mathfrak{sp}(2\times[7,9])$.
Pour cet exemple, $A'= A \setminus \{3\}$, et $f(A',D)$ est la colonne dont le double est
$$\aligned
dble(f(A',D))&=\begin{tabular}{|c|c|} \hline
\raisebox{-2pt}{$1$}&\raisebox{-2pt}{$4$}\\
\hline
\raisebox{-2pt}{$2$}&\raisebox{-2pt}{$5$}\\
\hline
\raisebox{-2pt}{$7$}&\raisebox{-2pt}{$7$}\\
\hline
\raisebox{-2pt}{$8$}&\raisebox{-2pt}{$9$}\\
\hline
\raisebox{-2pt}{$\overline{9}$}&\raisebox{-2pt}{$\overline 8$}\\
\hline
\raisebox{-2pt}{$\overline{5}$}&\raisebox{-2pt}{$\overline 3$}\\
\hline
\raisebox{-2pt}{$\overline{4}$}&\raisebox{-2pt}{$\overline 2$}\\
\hline
\raisebox{-2pt}{$\overline 3$}&\raisebox{-2pt}{$\overline 1$}\\
\hline
\end{tabular}\\
dble(f^{<6}(A^{<6},D^{<6}))&=\begin{tabular}{|c|c|} \hline
\raisebox{-2pt}{$1$}&\raisebox{-2pt}{$4$}\\
\hline
\raisebox{-2pt}{$2$}&\raisebox{-2pt}{$5$}\\
\hline
\raisebox{-2pt}{$\overline{5}$}&\raisebox{-2pt}{$\overline 3$}\\
\hline
\raisebox{-2pt}{$\overline{4}$}&\raisebox{-2pt}{$\overline 2$}\\
\hline
\raisebox{-2pt}{$\overline 3$}&\raisebox{-2pt}{$\overline 1$}\\
\hline
\end{tabular},\quad
dble(f^{>6}(A^{>6},D^{>6}))=\begin{tabular}{|c|c|}
\hline
\raisebox{-2pt}{$7$}&\raisebox{-2pt}{$7$}\\
\hline
\raisebox{-2pt}{$8$}&\raisebox{-2pt}{$9$}\\
\hline
\raisebox{-2pt}{$\overline{9}$}&\raisebox{-2pt}{$\overline 8$}\\
\hline
\end{tabular}.\endaligned
$$
L'entier $b_3=6$ n'est ni dans
$$
A^{<6}\cup B'^{<6}\cup C'^{<6} \cup D^{<6}=\{1,2\}\cup\{4,5\}\cup\{3,4,5\}\cup\{1,2,3\}
$$
ni dans
$$
A^{>6}\cup B^{>6}\cup C^{>6} \cup D^{>6}=\{7,8\}\cup\{7,9\}\cup\{9\}\cup\{8\}.
$$

\end{exe}

Le but de cet article est de montrer que l'ensemble des tableaux
quasi standards symplectiques forme une base de l'alg\`ebre de forme
r\'eduite qui respecte sa structure de $\mathfrak n^+$ module
ind\'ecomposable. Nous rappelons d'abord le jeu de taquin
symplectique d\'efini par J. T. Sheats dans \cite{SH}.


\subsection{Jeu de taquin symplectique}

\

Rappelons maintenant la d\'efinition du jeu de taquin symplectique de J. T. Sheats \cite{SH}.\\

Soit $T\setminus S$ un tableau de Young tordu de forme $\lambda\setminus\mu$. On d\'efinit le double de $T\setminus S$ en doublant les cases vides de $S$ et en doublant les bas remplis des colonnes comme ci-dessus. On dit que $T\setminus S$ est semi standard si $dble(T\setminus S)$ ainsi d\'efini est un tableau tordu semi standard. Voici un exemple :
$$
T\setminus S=\begin{tabular}{|c|c|c|}
\hline
 & \raisebox{-2pt}{$1$} & \raisebox{-2pt}{$2$}\\
\hline
 & \raisebox{-2pt}{$3$} & \raisebox{-2pt}{$4$}\\
\hline
 & \raisebox{-2pt}{$\overline{3}$} & \raisebox{-2pt}{$\overline{2}$}\\
\hline
\raisebox{-2pt}{$\overline{3}$} & \raisebox{-2pt}{$\overline{1}$}\\
\cline{1-2}
\raisebox{-2pt}{$\overline{2}$}\\
\cline{1-1}
\raisebox{-2pt}{$\overline{1}$}\\
\cline{1-1}
\end{tabular},\qquad dble(T\setminus S)=\begin{tabular}{|c|c|c|c|c|c|}
\hline
 & & \raisebox{-2pt}{$1$} & \raisebox{-2pt}{$2$} & \raisebox{-2pt}{$2$} & \raisebox{-2pt}{$3$}\\
\hline
 & & \raisebox{-2pt}{$3$} & \raisebox{-2pt}{$4$} & \raisebox{-2pt}{$4$} & \raisebox{-2pt}{$4$}\\
\hline
 & & \raisebox{-2pt}{$\overline{4}$} & \raisebox{-2pt}{$\overline{3}$} & \raisebox{-2pt}{$\overline{3}$} & \raisebox{-2pt}{$\overline{2}$}\\
\hline
\raisebox{-2pt}{$\overline{3}$} & \raisebox{-2pt}{$\overline{3}$} & \raisebox{-2pt}{$\overline{2}$} & \raisebox{-2pt}{$\overline{1}$}\\
\cline{1-4}
\raisebox{-2pt}{$\overline{2}$} & \raisebox{-2pt}{$\overline{2}$}\\
\cline{1-2}
\raisebox{-2pt}{$\overline{1}$} & \raisebox{-2pt}{$\overline{1}$}\\
\cline{1-2}
\end{tabular}
$$
Pour chaque colonne $c_j=f(A_j,D_j)=g(B_j,C_j)$ de $T\setminus S$, on note $\boxed{t_{ij}}$ la case $i$ de cette colonne et dans $dble(T\setminus S)$, la colonne $j$ devient deux colonnes. Les cases de ces colonnes sont not\'ees
$
\begin{tabular}{|c|c|}
\hline
\raisebox{-2pt}{$\alpha_{ij}$} & \raisebox{-2pt}{$\beta_{ij}$}\\
\hline
\end{tabular}$.

Lorsqu'on pointe un tableau semi standard tordu $T\setminus S$ en un coin int\'erieur de $S$, par convention on double la case $\boxed{\star}$ qui devient $\begin{array}{|c|c|}\hline \star&\star\\ \hline\end{array}$.\\

Le jeu de taquin symplectique consiste \`a partir d'un tableau semi standard tordu point\'e et \`a d\'eplacer la case point\'ee de la fa\c con suivante~: supposons que la case point\'ee soit en $(i,j)$. On note les parties remplies des colonnes de $T$ par $c_j=f(A_j,D_j)=g(B_j,C_j)$, alors

\begin{itemize}
\item[1] Si $(i,j+1)$ n'est pas une case de $T$ ou si $\beta_{(i+1)j}\leq \alpha_{i(j+1)}$, on permute la case point\'ee $\boxed{\star}$ de $T\setminus S$ avec la case $\boxed{t_{i+1,j}}$ imm\'ediatement en dessous, les autres cases restent inchang\'ees,\\
\item[2] Si $(i+1,j)$ n'est pas une case de $T$ ou si $\beta_{(i+1)j}> \alpha_{i(j+1)}$, on d\'eplace horizontalement la case point\'ee $\boxed{\star}$ suivant la r\`egle suivante~:
\subitem{\sl(i)} si $\alpha_{i,j+1}$ est non barr\'e, on remplace la colonne $c_j$ ainsi
$$
c_j=g(B_j,C_j)~~\longrightarrow~~c'_j=g(B_j\cup\{\alpha_{i,j+1}\},C_j)
$$
(la case point\'ee dispara\^\i t) et la colonne $c_{j+1}$ ainsi
$$
c_{j+1}=f(A_{j+1},D_{j+1})~~\longrightarrow~~c'_{j+1}=f(A_{j+1}\setminus\{\alpha_{i,j+1}\},D_{j+1})
$$
et la case $\boxed{\star}$ en $(i,j+1)$, les autres colonnes sont inchang\'ees.\\
\subitem{\sl(ii)} si $\alpha_{i,j+1}$ est barr\'e, on remplace la colonne $c_j$ ainsi
$$
c_j=f(A_j,D_j)~~\longrightarrow~~c'_j=f(A_j,D_j\cup\{\alpha_{i,j+1}\})
$$
(la case point\'ee dispara\^\i t) et la colonne $c_{j+1}$ ainsi
$$
c_{j+1}=g(B_{j+1},C_{j+1})~~\longrightarrow~~c'_{j+1}=g(B_{j+1},C_{j+1}\setminus\{\alpha_{i,j+1}\})
$$
et la case $\boxed{\star}$ en $(i,j+1)$, les autres colonnes sont inchang\'ees.\\
\item[3] Si ni $(i,j+1)$, ni $(i+1,j)$ n'est une case de $T$, le jeu s'arr\^ete.\\
\end{itemize}

\begin{exe}

Reprenons le tableau :
$$
T\setminus S=\begin{tabular}{|c|c|c|}
\hline
 & \raisebox{-2pt}{$1$} & \raisebox{-2pt}{$2$}\\
\hline
 & \raisebox{-2pt}{$3$} & \raisebox{-2pt}{$4$}\\
\hline
$\star$ & \raisebox{-2pt}{$\overline{3}$} & \raisebox{-2pt}{$\overline{2}$}\\
\hline
\raisebox{-2pt}{$\overline{3}$} & \raisebox{-2pt}{$\overline{1}$}\\
\cline{1-2}
\raisebox{-2pt}{$\overline{2}$}\\
\cline{1-1}
\raisebox{-2pt}{$\overline{1}$}\\
\cline{1-1}
\end{tabular}
$$
Le jeu de taquin donne successivement :
$$\aligned
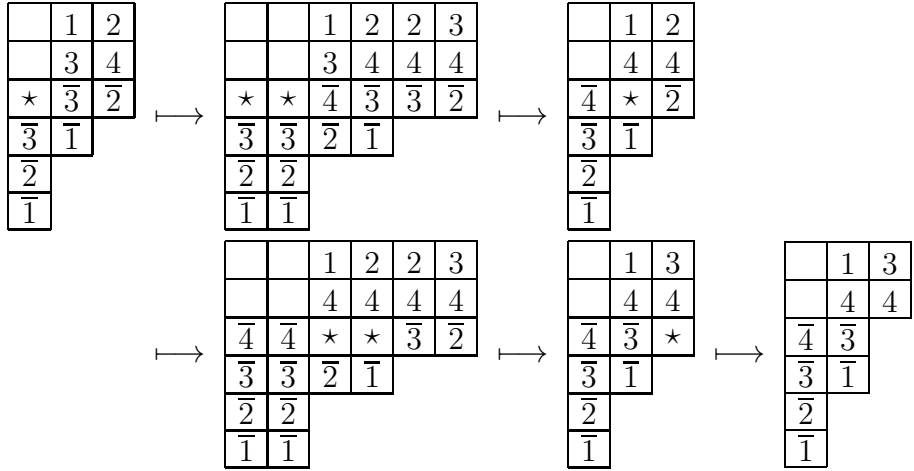
\begin{tabular}{|c|c|c|}
\hline
 & \raisebox{-2pt}{$1$} & \raisebox{-2pt}{$2$}\\
\hline
 & \raisebox{-2pt}{$3$} & \raisebox{-2pt}{$4$}\\
\hline
$\star$ & \raisebox{-2pt}{$\overline{3}$} & \raisebox{-2pt}{$\overline{2}$}\\
\hline
\raisebox{-2pt}{$\overline{3}$} & \raisebox{-2pt}{$\overline{1}$}\\
\cline{1-2}
\raisebox{-2pt}{$\overline{2}$}\\
\cline{1-1}
\raisebox{-2pt}{$\overline{1}$}\\
\cline{1-1}
\end{tabular}&~\longmapsto~
\begin{tabular}{|c|c|c|c|c|c|}
\hline
 & & \raisebox{-2pt}{$1$} & \raisebox{-2pt}{$2$} & \raisebox{-2pt}{$2$} & \raisebox{-2pt}{$3$}\\
\hline
 & & \raisebox{-2pt}{$3$} & \raisebox{-2pt}{$4$} & \raisebox{-2pt}{$4$} & \raisebox{-2pt}{$4$}\\
\hline
$\star$ & $\star$ & \raisebox{-2pt}{$\overline{4}$} & \raisebox{-2pt}{$\overline{3}$} & \raisebox{-2pt}{$\overline{3}$} & \raisebox{-2pt}{$\overline{2}$}\\
\hline
\raisebox{-2pt}{$\overline{3}$} & \raisebox{-2pt}{$\overline{3}$} & \raisebox{-2pt}{$\overline{2}$} & \raisebox{-2pt}{$\overline{1}$}\\
\cline{1-4}
\raisebox{-2pt}{$\overline{2}$} & \raisebox{-2pt}{$\overline{2}$}\\
\cline{1-2}
\raisebox{-2pt}{$\overline{1}$} & \raisebox{-2pt}{$\overline{1}$}\\
\cline{1-2}
\end{tabular}~\longmapsto~\begin{tabular}{|c|c|c|}
\hline
 & \raisebox{-2pt}{$1$} & \raisebox{-2pt}{$2$}\\
\hline
 & \raisebox{-2pt}{$4$} & \raisebox{-2pt}{$4$}\\
\hline
\raisebox{-2pt}{$\overline{4}$} & $\star$ & \raisebox{-2pt}{$\overline{2}$}\\
\hline
\raisebox{-2pt}{$\overline{3}$} & \raisebox{-2pt}{$\overline{1}$}\\
\cline{1-2}
\raisebox{-2pt}{$\overline{2}$}\\
\cline{1-1}
\raisebox{-2pt}{$\overline{1}$}\\
\cline{1-1}
\end{tabular}\\
&~\longmapsto~
\begin{tabular}{|c|c|c|c|c|c|}
\hline
 & & \raisebox{-2pt}{$1$} & \raisebox{-2pt}{$2$} & \raisebox{-2pt}{$2$} & \raisebox{-2pt}{$3$}\\
\hline
 & & \raisebox{-2pt}{$4$} & \raisebox{-2pt}{$4$} & \raisebox{-2pt}{$4$} & \raisebox{-2pt}{$4$}\\
\hline
\raisebox{-2pt}{$\overline{4}$} & \raisebox{-2pt}{$\overline{4}$} & $\star$ & $\star$ & \raisebox{-2pt}{$\overline{3}$} & \raisebox{-2pt}{$\overline{2}$}\\
\hline
\raisebox{-2pt}{$\overline{3}$} & \raisebox{-2pt}{$\overline{3}$} & \raisebox{-2pt}{$\overline{2}$} & \raisebox{-2pt}{$\overline{1}$}\\
\cline{1-4}
\raisebox{-2pt}{$\overline{2}$} & \raisebox{-2pt}{$\overline{2}$}\\
\cline{1-2}
\raisebox{-2pt}{$\overline{1}$} & \raisebox{-2pt}{$\overline{1}$}\\
\cline{1-2}
\end{tabular}~\longmapsto~\begin{tabular}{|c|c|c|}
\hline
 & \raisebox{-2pt}{$1$} & \raisebox{-2pt}{$3$}\\
\hline
 & \raisebox{-2pt}{$4$} & \raisebox{-2pt}{$4$}\\
\hline
\raisebox{-2pt}{$\overline{4}$} & \raisebox{-2pt}{$\overline{3}$} & $\star$\\
\hline
\raisebox{-2pt}{$\overline{3}$} & \raisebox{-2pt}{$\overline{1}$}\\
\cline{1-2}
\raisebox{-2pt}{$\overline{2}$}\\
\cline{1-1}
\raisebox{-2pt}{$\overline{1}$}\\
\cline{1-1}
\end{tabular}
~\longmapsto~\begin{tabular}{|c|c|c|}
\hline
 & \raisebox{-2pt}{$1$} & \raisebox{-2pt}{$3$}\\
\hline
 & \raisebox{-2pt}{$4$} & \raisebox{-2pt}{$4$}\\
\hline
\raisebox{-2pt}{$\overline{4}$} & \raisebox{-2pt}{$\overline{3}$}\\
\cline{1-2}
\raisebox{-2pt}{$\overline{3}$} & \raisebox{-2pt}{$\overline{1}$}\\
\cline{1-2}
\raisebox{-2pt}{$\overline{2}$}\\
\cline{1-1}
\raisebox{-2pt}{$\overline{1}$}\\
\cline{1-1}
\end{tabular}\\
\endaligned
$$

Ensuite, on peut recommencer avec le tableau obtenu :
$$
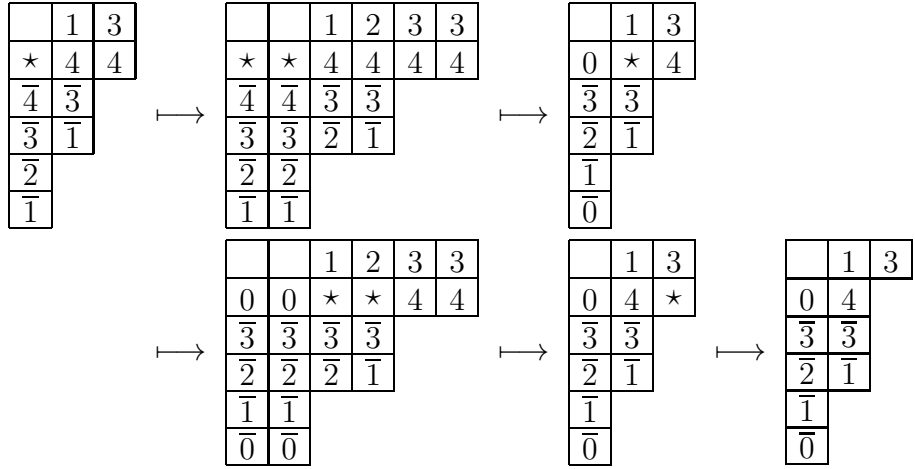
\aligned
\begin{tabular}{|c|c|c|}
\hline
 & \raisebox{-2pt}{$1$} & \raisebox{-2pt}{$3$}\\
\hline

$\star$ & \raisebox{-2pt}{$4$} & \raisebox{-2pt}{$4$}\\
\hline
\raisebox{-2pt}{$\overline{4}$} & \raisebox{-2pt}{$\overline{3}$}\\
\cline{1-2}
\raisebox{-2pt}{$\overline{3}$} & \raisebox{-2pt}{$\overline{1}$}\\
\cline{1-2}
\raisebox{-2pt}{$\overline{2}$}\\
\cline{1-1}
\raisebox{-2pt}{$\overline{1}$}\\
\cline{1-1}
\end{tabular}&~\longmapsto~\begin{tabular}{|c|c|c|c|c|c|}
\hline
 & & \raisebox{-2pt}{$1$} & \raisebox{-2pt}{$2$} & \raisebox{-2pt}{$3$} & \raisebox{-2pt}{$3$}\\
\hline
$\star$ & $\star$ & \raisebox{-2pt}{$4$} & \raisebox{-2pt}{$4$} & \raisebox{-2pt}{$4$} & \raisebox{-2pt}{$4$}\\
\hline
\raisebox{-2pt}{$\overline{4}$} & \raisebox{-2pt}{$\overline{4}$} & \raisebox{-2pt}{$\overline{3}$} & \raisebox{-2pt}{$\overline{3}$}\\
\cline{1-4}
\raisebox{-2pt}{$\overline{3}$} & \raisebox{-2pt}{$\overline{3}$} & \raisebox{-2pt}{$\overline{2}$} & \raisebox{-2pt}{$\overline{1}$}\\
\cline{1-4}
\raisebox{-2pt}{$\overline{2}$} & \raisebox{-2pt}{$\overline{2}$}\\
\cline{1-2}
\raisebox{-2pt}{$\overline{1}$} & \raisebox{-2pt}{$\overline{1}$}\\
\cline{1-2}
\end{tabular}~\longmapsto~\begin{tabular}{|c|c|c|}
\hline
 & \raisebox{-2pt}{$1$} & \raisebox{-2pt}{$3$}\\
\hline
\raisebox{-2pt}{$0$} & $\star$ & \raisebox{-2pt}{$4$}\\
\hline
\raisebox{-2pt}{$\overline{3}$} & \raisebox{-2pt}{$\overline{3}$}\\
\cline{1-2}
\raisebox{-2pt}{$\overline{2}$} & \raisebox{-2pt}{$\overline{1}$}\\
\cline{1-2}
\raisebox{-2pt}{$\overline{1}$}\\
\cline{1-1}
\raisebox{-2pt}{$\overline{0}$}\\
\cline{1-1}
\end{tabular}\\
&~\longmapsto~\begin{tabular}{|c|c|c|c|c|c|}
\hline
 & & \raisebox{-2pt}{$1$} & \raisebox{-2pt}{$2$} & \raisebox{-2pt}{$3$} & \raisebox{-2pt}{$3$}\\
\hline
\raisebox{-2pt}{$0$} & \raisebox{-2pt}{$0$} & $\star$ & $\star$ & \raisebox{-2pt}{$4$} & \raisebox{-2pt}{$4$}\\
\hline
\raisebox{-2pt}{$\overline{3}$} & \raisebox{-2pt}{$\overline{3}$} & \raisebox{-2pt}{$\overline{3}$} & \raisebox{-2pt}{$\overline{3}$}\\
\cline{1-4}
\raisebox{-2pt}{$\overline{2}$} & \raisebox{-2pt}{$\overline{2}$} & \raisebox{-2pt}{$\overline{2}$} & \raisebox{-2pt}{$\overline{1}$}\\
\cline{1-4}
\raisebox{-2pt}{$\overline{1}$} & \raisebox{-2pt}{$\overline{1}$}\\
\cline{1-2}
\raisebox{-2pt}{$\overline{0}$} & \raisebox{-2pt}{$\overline{0}$}\\
\cline{1-2}
\end{tabular}~\longmapsto~\begin{tabular}{|c|c|c|}
\hline
 & \raisebox{-2pt}{$1$} & \raisebox{-2pt}{$3$}\\
\hline
\raisebox{-2pt}{$0$} & \raisebox{-2pt}{$4$} & $\star$\\
\hline
\raisebox{-2pt}{$\overline{3}$} & \raisebox{-2pt}{$\overline{3}$}\\
\cline{1-2}
\raisebox{-2pt}{$\overline{2}$} & \raisebox{-2pt}{$\overline{1}$}\\
\cline{1-2}
\raisebox{-2pt}{$\overline{1}$}\\
\cline{1-1}
\raisebox{-2pt}{$\overline{0}$}\\
\cline{1-1}
\end{tabular}
~\longmapsto~\begin{tabular}{|c|c|c|}
\hline
 & \raisebox{-2pt}{$1$} & \raisebox{-2pt}{$3$}\\
\hline
\raisebox{-2pt}{$0$} & \raisebox{-2pt}{$4$}\\
\cline{1-2}
\raisebox{-2pt}{$\overline{3}$} & \raisebox{-2pt}{$\overline{3}$}\\
\cline{1-2}
\raisebox{-2pt}{$\overline{2}$} & \raisebox{-2pt}{$\overline{1}$}\\
\cline{1-2}
\raisebox{-2pt}{$\overline{1}$}\\
\cline{1-1}
\raisebox{-2pt}{$\overline{0}$}\\
\cline{1-1}
\end{tabular}\\
\endaligned
$$
\end{exe}

Dans cet exemple, on voit que le nombre 0 peut appara\^\i tre (cf. \cite{SH}). En fait J. T. Sheats a montr\'e qu'il ne peut appara\^\i tre que dans la premi\`ere colonne, et qu'elle appara\^\i t en m\^eme temps que $\overline{0}$.\\


\subsection{Jeu de taquin et tableaux non quasi standards}

\

Appliquons le jeu de taquin \`a un tableau semi standard symplectique $T$ qui n'est pas quasi standard symplectique en $s$, c'est \`a dire que $T$ appartient \`a $NQS_s$ et poss\`ede une colonne de hauteur $s$.\\

On ajoute \`a gauche de $T$ une colonne triviale $c_0$ de hauteur
$n$ dont on vide les $s$ premi\`eres cases (on note $T_0\setminus S$
le tableau obtenu). On pointe le coin inf\'erieur de $S$ et on
applique le jeu de taquin.

\begin{prop}
Lorsqu'on applique le jeu de taquin symplectique \`a $T_0\setminus
S$, les \'etoiles se d\'eplacent toujours horizontalement de la
gauche vers la droite, l'indice 0 n'appara\^\i t pas et le tableau
obtenu a pour premi\`ere colonne la colonne triviale $c_0$ \`a qui
on a vid\'e les $s-1$ premi\`eres cases. Si $s>1$, le tableau $T'$
form\'e par les colonnes suivantes est semi standard, non quasi
standard en $s-1$ et poss\`ede une colonne de hauteur $s-1$.
\end{prop}

\noindent
{\bf Preuve}

Par construction, la colonne $c_0\setminus S$ se double en $(c_0\setminus S)(c_0\setminus S)$, \`a droite de $\boxed{\star}$, il y a $\boxed{s}$ ($s$ n'est pas barr\'e) et au dessous $\boxed{s+1}$. Le premier pas du jeu de taquin consiste simplement \`a permuter les cases $\boxed{\star}$ et $\boxed{s}$ des colonnes 0 et 1. En particulier, 0 n'appara\^\i t pas et la premi\`ere colonne a la forme annonc\'ee. La colonne $c_1=f(A_1,D_1)=g(B_1,C_1)$ devient $c'_1=f(A'_1,D_1)$, point\'ee en $s$.\\

Supposons qu'apr\`es un certain nombre de pas, la case point\'ee soit toujours sur la ligne $s$, dans la colonne $i+1$, on fait les hypoth\`eses suivantes :
\begin{itemize}
\item[(H1)] Les $i$ premi\`eres colonnes de notre nouveau tableau, not\'ees $c"_1$, \dots, $c"_i$ sont de la forme $f(A"_j,D"_j)=g(B"_j,C"_j)$, le tableau $c"_1\dots c"_i$ est dans $NQS_{s-1}$.
\item[(H2)] La colonne $i+1$ est devenue $c'_{i+1}$, elle contient une \'etoile \`a la ligne $s$, dans $dble(T)$, on a $t_{s-1,2(i+1)}<t_{s,2i+1}$.
\end{itemize}
Les colonnes suivantes $i+2,\dots$ n'ont pas \'et\'e modifi\'ees, on les note $c_j$ ($j>i+1$). On repr\'esente cette situation par le sh\'ema suivant :
$$
\begin{array}{|c|c|c|}
(c"_i)&(c'_{i+1})&(c_{i+2})\\
&&\\
&&\\
\cline{1-2}
\vspace{-0.2cm}&\vspace{-0.2cm}&\vspace{-0.2cm}\\
\cline{1-2}
&\star&\\
\cline{2-3}
\vspace{-0.2cm}&\vspace{-0.2cm}&\vspace{-0.2cm}\\
\cline{2-3}
&&\\
&&\\
\end{array}
$$

\underline{Les \'etoiles se d\'eplacent vers la droite}\\

Si $s\leq\#A_{i+2}$, dans ce cas, on a, gr\^ace au lemme 2, la situation suivante~:
$$
dble=\begin{array}{|c|c|}
dble(c'_{i+1})&dble(c_{i+2})\\
&\\
a_{s-1,i+1}~\hspace{0.5cm}~b'_{s-1,i+1}&\\
&\\
\cline{1-1}
\vspace{-0.2cm}&\vspace{-0.2cm}\\
\cline{1-1}
&\\
\star\quad\quad\star&a_{s,i+2}\hspace{2cm}\\
&\\
\cline{1-2}
\vspace{-0.2cm}&\vspace{-0.2cm}\\
\cline{1-2}
&\\
~\hfill b_{s+1,i+1}&\\
\hfil \text{ou}&\\
\hfill \overline{c}_{t-1,i+1}&\\
&\\
\end{array}
$$
Le d\'eplacement suivant est horizontal.\\

Si $s>\#A_{i+2}$, dans ce cas, on a, gr\^ace au lemme 2 ou \`a la remarque qui le suit, la situation suivante~:
$$
dble=\begin{array}{|c|c|}
dble(c'_{i+1})&dble(c_{i+2})\\
&\\
&\\
\cline{1-1}
\vspace{-0.2cm}&\vspace{-0.2cm}\\
\cline{1-1}
&\\
\star\quad\quad\star&\overline{c}_{t,i+2}\hspace{2cm}\\
&\\
\cline{1-2}
\vspace{-0.2cm}&\vspace{-0.2cm}\\
\cline{1-2}
&\\
~\hfill \overline{d}_{t-1,i+1}&\\
&\\
\end{array}
$$
Le d\'eplacement suivant est horizontal.\\

\underline{Au pas suivant le tableau form\'e des colonnes $0,\dots,i+1$ est dans $NQS_{s-1}$}\\

Si $s\leq\# A_{i+2}$, dans ce cas, on a, gr\^ace au lemme 3, la situation suivante~:
$$
dble=\begin{array}{|c|c|}
dble(c"_i)&dble(c'_{i+1})\\
&\\
&a_{s-1,i+1}\quad\quad b'_{s-1,i+1}(<b_{s,i+1})\\
&\\
\cline{1-2}
\vspace{-0.2cm}&\vspace{-0.2cm}\\
\cline{1-2}
&\\
\qquad a_{s,i+1}&(b_{s,i+1}\leq)a"_{s,i+1}\quad b"_{s,i+1}(=a_{s,i+2})\\
&\\
\cline{2-2}
\vspace{-0.2cm}&\vspace{-0.2cm}\\
\cline{2-2}
\end{array}
$$

\vskip0.5cm
Ce qui nous donne les deux in\'egalit\'es demand\'ees $a"_{s-1,i+1}=a_{s-1,i+1}<b"_{s,i}=a_{s,i+1}$ et $b"_{s-1,i+1}=b_{s-1,i+1}<b_{s,i+1}\leq a"_{s,i+1}$.\\

Si $\# A_{i+2}<s\leq\#A_{i+1}$, dans ce cas, on a, gr\^ace au lemme 3 et \`a la remarque qui le suit, la situation suivante~:
$$
dble=\begin{array}{|c|c|}
dble(c"_i)&dble(c'_{i+1})\\
&\\
&a_{s-1,i+1}\quad\quad b'_{s-1,i+1}(<b_{s,i+1})\\
&\\
\cline{1-2}
\vspace{-0.2cm}&\vspace{-0.2cm}\\
\cline{1-2}
&\\
\qquad a_{s,i+1}&\overline{c}"_{t,i+1}\quad \overline{d}"_{t,i+1}(=\overline{c}_{t,i+2})\\
&\\
\cline{2-2}
\vspace{-0.2cm}&\vspace{-0.2cm}\\
\cline{2-2}
\end{array}
$$
\vskip0.5cm
Ce qui nous donne les deux in\'egalit\'es demand\'ees $a"_{s-1,i+1}=a_{s-1,i+1}<b"_{s,i}=a_{s,i+1}$ et $b"_{s-1,i+1}<\overline{c}"_{t,i+1}$.\\

Si $\# A_{i+1}<s$, dans ce cas, on a, gr\^ace au lemme 3 et \`a la remarque qui le suit, les deux situations suivantes~:
$$
dble=\begin{array}{|c|c|}
dble(c"_i)&dble(c'_{i+1})\\
&\\
&a_{s-1,i+1}\quad\quad b'_{s-1,i+1}(<b_{s,i+1})\\
&\\
\cline{1-2}
\vspace{-0.2cm}&\vspace{-0.2cm}\\
\cline{1-2}
&\\
\qquad \overline{c}_{t,i+1}&\overline{c}"_{t,i+1}\quad \overline{d}"_{t,i+1}(=\overline{c}_{t,i+2})\\
&\\
\cline{2-2}
\vspace{-0.2cm}&\vspace{-0.2cm}\\
\cline{2-2}
\end{array}
$$
\vskip0.5cm
Ce qui nous donne les deux in\'egalit\'es demand\'ees $a"_{s-1,i+1}<\overline{c}"_{t,i}$ et $b"_{s-1,i+1}<\overline{c}"_{t,i+1}$.\\

Ou bien
$$
dble=\begin{array}{|c|c|}
dble(c"_i)&dble(c'_{i+1})\\
&\\
&\overline{c}"_{t+1,i+1}\quad\quad \overline{d}"_{t+1,i+1}(<\overline{d}_{t,i+1})\\
&\\
\cline{1-2}
\vspace{-0.2cm}&\vspace{-0.2cm}\\
\cline{1-2}
&\\
\qquad \overline{c}_{t,i+1}&(\overline{d}_{t,i+1}\leq)\overline{c}"_{t,i+1}\quad \overline{d}"_{t,i+1}(=\overline{c}_{t,i+2})\\
&\\
\cline{2-2}
\vspace{-0.2cm}&\vspace{-0.2cm}\\
\cline{2-2}
\end{array}
$$
Ce qui nous donne les deux in\'egalit\'es demand\'ees
$\overline{c}"_{t+1,i+1}\leq\overline{d}"_{t+1,i+1}<\overline{d}_{t,i+1}\leq\overline{c}_{t,i+1}=\overline{d}"_{t,i}$
et
$\overline{d}"_{t+1,i+1}<\overline{d}_{t,i+1}\leq\overline{c}"_{t,i+1}$.

$\hfill\square$

Soit $T$ un tableau de $SS^{<\lambda>}$ qui n'est pas quasi
standard. Soit $s$ un entier tel que $T\in NQS_s$ et $T$ poss\`ede
une colonne de hauteur $s$. On notera ceci : $T\in
NQS_s^{<\lambda>}$. Supposons $T\notin NQS_t^{<\lambda>}$, pour tout $t>s$.
On ajoute \`a $T$ une colonne triviale avec $s$ cases vides
$c_0\setminus S$, on applique le jeu de taquin symplectique, on
retire la premi\`ere colonne (triviale avec $s-1$ cases vides) et on
obtient un tableau $T'\in NQS_{s-1}^{<\lambda-[s]+[s-1]>}$ o\`u
$[s]$ d\'esigne le $n$-uplet $(0,\dots,1,\dots,0)$, le 1 \'etant \`a
la $s^{i\grave{e}me}$ place. On notera $T'=sjdt_s(T)$. Il est
possible que $T'$ soit dans $NQS_s^{<\lambda-[s]+[s-1]>}$. Cependant
$T'$ ne peut pas \^etre dans $NQS_t^{<\lambda-[s]+[s-1]>}$, avec
$t>s$.

\begin{lem}
Pour tout $t>s$, si $T$ n'est pas dans $NQS_t$, alors $T'=sjdt_s(T)$ n'est pas non plus dans $NQS_t$.
\end{lem}

\noindent
{\bf Preuve}

En effet, si $T$ n'a pas de colonne de hauteur $t$, $T'$ n'en n'a pas non plus. Si $T$ a une colonne de hauteur $t$ et n'est pas dans $NQS_t$, le double de $T$ est tel que pour chaque $t$ il existe un `blocage' de la forme $t_{t,j}\geq t_{t+1,j-1}$. Montrons que ce blocage ne dispara\^\i t pas au cours du jeu de taquin symplectique. Supposons qu'\`a la colonne $c_i$ de $T$, on ait $s\leq\#A$. Si on ajoute une entr\'ee barr\'ee $\overline{u}$, on a vu que la partie de la nouvelle colonne $dble(c"_i)$ situ\'ee en dessous de la ligne $s$ est le double de $f(A^{<b_s},D^{<b_s})$, c'est \`a dire co\"\i ncide avec la nouvelle colonne $dble(c_i)$ situ\'ee en dessous de la ligne $s$. S'il y avait un blocage, il n'a pas disparu. Si on ajoute une entr\'ee $u$ qui n'est pas barr\'ee, on est dans la situation suivante (la parenth\`ese signifie que la ligne correspondante peut exister ou ne pas exister), si $u=b_s$, on obtient $c"_i=c_i$ et aucun blocage ne dispara\^\i t, si $u>b_s$ et $u\notin C$,
$$\aligned
&dbl(c_i)=\begin{array}{cc}
A^{<b_s}&B^{<b_s}\\
a_s&b_s\\
A^{>b_s}&B^{>b_s}\\
\overline{C}^{>b_s}&\overline{D}^{>b_s}\\
\left(\overline{b_s}\right.&\left.\overline{a_s}\right)\\
\overline{C}^{<b_s}&\overline{D}^{<b_s}
\end{array}~~\mapsto~~dble(c'_i)=\begin{array}{cc}
A^{<b_s}&B^{<b_s}\\
\star&\star\\
A^{>b_s}&B^{>b_s}\\
\overline{C}^{>b_s}&\overline{D}^{>b_s}\\
\left(\overline{a_s}\right.&\left.\overline{a_s}\right)\\
\overline{C}^{<b_s}&\overline{D}^{<b_s}
\end{array}~~\mapsto\\
&\mapsto~~dble(c"_i)=\begin{array}{cc}
A^{<b_s}&B^{<b_s}\\
u&u\\
A^{>b_s}&B^{>b_s}\\
\overline{C}^{>b_s}&\overline{D}^{>b_s}\\
\left(\overline{a_s}\right.&\left.\overline{a_s}\right)\\
\overline{C}^{<b_s}&\overline{D}^{<b_s}
\end{array}
\endaligned
$$
le seul changement \'eventuel, en dessous de la ligne $s$, est le remplacement de $\overline{b_s}$ par $\overline{a_s}$. Ce remplacement n'appporte aucun nouveau d\'eblocage. Enfin si $u>b_s$ appartient \`a $C$, on a $u=c_a$, par construction $u$ est le plus petit \'el\'ement de $J_i"^{\geq b_s}$, on a, avec nos notations, $D"^{\geq b_s}=D'^{\geq b_s}\setminus\{u\}\cup\{v=d"_b\}$, le seul changement des colonnes en dessous de $s$, \`a part le changement \'eventuel de $\overline{b_s}$ en $\overline{a_s}$, est la partie comprise entre les ligne barr\'ees d'indices $a$ et $b$. Plus pr\'ecis\'ement, cette partie devient :
$$
\begin{array}{ccc}
\overline{c}_{a+1}&&\overline{d}_{a+1}\\
\overline{c}_{a}&&\overline{c}_{a}\\
\vdots&&\vdots\\
\overline{c}_{b+1}&&\overline{c}_{b+1}\\
\overline{c}_{b}&&\overline{c}_{b}\\
\end{array}~~\mapsto~~
\begin{array}{ccc}
\overline{c}_{a+1}&&\overline{d}_{a+1}\\
\overline{c}_{a}&&\overline{c}_{a-1}\\
\vdots&&\vdots\\
\overline{c}_{b+1}&&\overline{c}_{b}\\
\overline{c}_{b}&&\overline{d"}_{b}\\
\end{array}
$$
On voit appara\^\i tre des blocages entre ces deux colonnes entre les lignes $a$ et $b$. Aucun ancien blocage ne dispara\^\i t. Le m\^eme argument s'applique si $s>\#A$.

$\hfill\square$

On peut maintenant r\'ep\'eter le jeu de taquin symplectique sur
$T'=sjdt_s(T)$. Si $T'$ n'est pas quasi standard (c'est en
particulier le cas si $s>1$), il existe $s'\leq s$ tel que $T'\in
NQS_{s'}^{<\lambda\setminus[s]\cup[s-1]>}$ et $T'\notin NQS_{t'}^{<\lambda\setminus[s]\cup[s-1]>}$, pour tout $t'>s'$, on construit
$T"=sjdt_{s'}(T')$, etc\dots Au bout d'un nombre fini
d'op\'erations, on obtient un tableau $\varphi(T)$ quasi standard :
$\varphi(T)\in QS^{<\mu>}$, et on d\'efinit ainsi une application
$\varphi$ de $SS^{<\lambda>}$ dans $\sqcup_{\mu\subset \lambda}
QS^{<\mu>}$.

\begin{theo}

\

L'application $\varphi$ est bijective de $SS^{<\lambda>}$ sur $\sqcup_{\mu\subset \lambda} QS^{<\mu>}$.
\end{theo}

\noindent
{\bf Preuve}

D'apr\`es le th\'eor\`eme 7.3 de \cite{SH}, on sait que le jeu de taquin symplectique $sjdt$ est injectif et que son application
inverse est de la forme $\sigma\circ sjdt\circ\sigma$ o\`u $\sigma$ est le retournement d'un tableau, accompagn\'e du changement des
entr\'ees barr\'ees en non barr\'ees et des non barr\'ees en barr\'ees. L'application $sjdt_s$ d\'eplace l'\'etoile vers la droite jusqu'\`a la derni\`ere case de la ligne $s$. On r\'ep\`ete cette op\'eration pour r\'ealiser $\varphi$. On obtient un tableau $\varphi(T)$ de forme $\mu$ et des \'etoiles succesives \`a droite de ce tableau qui remplissent le tableau tordu de forme $\lambda\setminus \mu$ de bas en haut et de droite \`a gauche (on remplit les lignes successivement en commen\c cant par la derni\`ere et dans chaque ligne de droite \`a gauche).

Si maintenant $T$ est un tableau quelconque de $SS^{<\mu>}$ avec
$\mu\subset\lambda$, on lui ajoute \`a gauche autant de colonnes
triviales qu'il y a de cases sur la premi\`ere ligne de
$\lambda\setminus\mu$ (disons $d$ colonnes), et le tableau tordu de
forme $\lambda\setminus \mu$ en haut \`a droite, on remplit ce
tableau tordu par des \'etoiles num\'erot\'ees comme ci-dessus, et
on applique le jeu de taquin inverse. On obtient par construction un
tableau $\theta(T)$ de forme $(\lambda\cup
d[n])\setminus(\lambda\setminus \mu)$, puisque d'apr\`es le
th\'eor\`eme 7.3 de \cite{SH}, les chemins successifs des \'etoiles
ne se croisent pas (au sens de \cite{SH}) : les derni\`eres
\'etoiles sont sur la premi\`ere ligne de notre tableau, le jeu de
taquin inverse les ram\`ene succesivement, dans l'ordre d\'ecroissant
le long de la premi\`ere ligne, le plus \`a gauche possible. Les \'etoiles suivantes sont sur la
ligne 2. Elles ne peuvent pas passer par la ligne 1, puisque les
chemins ne se croisent pas. Elles reviennent donc, le plus \`a
gauche possible, le long de cette ligne, etc\dots  A cause
de la forme de notre tableau, 0 n'appara\^\i t
jamais. En effet pour que 0 apparaisse, il faut qu'\`a un moment
donn\'e il y ait $\overline{1}$ \`a gauche de l'\'etoile et 1 au
dessus de l'\'etoile, donc il y a une case \`a gauche de ce 1 qui
contient n\'ecessairement 1. Mais alors la colonne \`a gauche de
l'\'etoile et celle au dessus de l'\'etoile forment un tableau qui
n'est pas semi standard, ceci est impossible d'apr\`es \cite{SH}. Le
tableau $\theta(T)$ obtenu est donc semi standard et par
construction ses $d$ premi\`eres colonnes sont des bas de colonnes
triviales. Ensuite, on compl\`ete le tableau $\theta(T)$ en
compl\`etant les $d$ premi\`eres colonnes en des colonnes triviales.
On supprime les $d$ premi\`eres colonnes triviales et on obtient un
tableau $\psi(T)$ semi standard de forme $\lambda$. D'apr\`es
\cite{SH}, $\varphi(\psi(T))=T$, $\varphi$ est bijective.

$\hfill\square$
\begin{exe} Cas de $\mathfrak{sp}(8)$
$$
T=\begin{tabular}{|c|c|c|} \hline
\raisebox{-2pt}{$1$} & \raisebox{-2pt}{$1$} & \raisebox{-2pt}{$3$}\\
\hline
\raisebox{-2pt}{$2$} & \raisebox{-2pt}{$3$} & \raisebox{-2pt}{$\overline{3}$}\\
\hline
\raisebox{-2pt}{$3$} & \raisebox{-2pt}{$\overline{3}$}\\
\cline{1-2}
\raisebox{-2pt}{$\overline{3}$}\\
\cline{1-1}
\end{tabular}\\~~ \Longrightarrow ~~\hbox{double}~(T)=\begin{tabular}{|c|c|c|c|c|c|} \hline
\raisebox{-2pt}{$1$} & \raisebox{-2pt}{$1$} & \raisebox{-2pt}{$1$}& \raisebox{-2pt}{$1$}& \raisebox{-2pt}{$3$}& \raisebox{-2pt}{$4$}\\
\hline
\raisebox{-2pt}{$2$} &\raisebox{-2pt}{$2$}& \raisebox{-2pt}{$3$}& \raisebox{-2pt}{$4$}& \raisebox{-2pt}{$\overline{4}$} & \raisebox{-2pt}{$\overline{3}$}\\
\hline
\raisebox{-2pt}{$3$} &\raisebox{-2pt}{$4$}&\raisebox{-2pt}{$\overline{4}$}& \raisebox{-2pt}{$\overline{3}$}\\
\cline{1-4}
\raisebox{-2pt}{$\overline{4}$}&\raisebox{-2pt}{$\overline{3}$}\\
\cline{1-2}
\end{tabular}.$$
Puisque $4<\overline{4}<\overline{3}$, le tableau $T\in
NQS_3^{(0,1,1,1)}$ . En oubliant l'ajout initial et le retrait final
des colonnes triviales, on a successivement:
$$\begin{tabular}{|c|c|c|} \hline
 & \raisebox{-2pt}{$1$} & \raisebox{-2pt}{$3$}\\
\hline
 & \raisebox{-2pt}{$3$} & \raisebox{-2pt}{$\overline{3}$}\\
\hline
$\star$ & \raisebox{-2pt}{$\overline{3}$}\\
\cline{1-2}
\raisebox{-2pt}{$\overline{3}$}\\
\cline{1-1}
\end{tabular}~~\longmapsto \begin{tabular}{|c|c|c|} \hline
 & \raisebox{-2pt}{$1$} & \raisebox{-2pt}{$3$}\\
\hline
 & \raisebox{-2pt}{$4$} & \raisebox{-2pt}{$\overline{3}$}\\
\hline
\raisebox{-2pt}{$\overline{4}$}\\
\cline{1-1}
\raisebox{-2pt}{$\overline{3}$}\\
\cline{1-1}
\end{tabular}~~
\longmapsto\begin{tabular}{|c|c|c|} \hline
 & \raisebox{-2pt}{$1$} & \raisebox{-2pt}{$3$}\\
\hline
 $\star$& \raisebox{-2pt}{$4$} & \raisebox{-2pt}{$\overline{3}$}\\
\hline
\raisebox{-2pt}{$\overline{4}$}\\
\cline{1-1}
\raisebox{-2pt}{$\overline{3}$}\\
\cline{1-1}
\end{tabular}~~\longmapsto
\begin{tabular}{|c|c|c|} \hline
 & \raisebox{-2pt}{$1$} & \raisebox{-2pt}{$3$}\\
\hline
 \raisebox{-2pt}{$2$}&$\star$ & \raisebox{-2pt}{$\overline{3}$}\\
\hline
\raisebox{-2pt}{$\overline{3}$}\\
\cline{1-1}
\raisebox{-2pt}{$\overline{2}$}\\
\cline{1-1}
\end{tabular}
~~\longmapsto \begin{tabular}{|c|c|c|} \hline
 & \raisebox{-2pt}{$1$} & \raisebox{-2pt}{$4$}\\
\hline
 \raisebox{-2pt}{$2$}& \raisebox{-2pt}{$\overline{4}$}\\
\cline{1-2}
\raisebox{-2pt}{$\overline{3}$}\\
\cline{1-1}
\raisebox{-2pt}{$\overline{2}$}\\
\cline{1-1}
\end{tabular}~~$$
$$\longmapsto \begin{tabular}{|c|c|c|} \hline
& \raisebox{-2pt}{$1$} & \raisebox{-2pt}{$4$}\\
\hline
 $\star$& \raisebox{-2pt}{$\overline{4}$}\\
\cline{1-2}
\raisebox{-2pt}{$\overline{3}$}\\
\cline{1-1}
\raisebox{-2pt}{$\overline{2}$}\\
\cline{1-1}
\end{tabular}~~\longmapsto
\begin{tabular}{|c|c|c|} \hline
& \raisebox{-2pt}{$1$} & \raisebox{-2pt}{$4$}\\
\hline
 \raisebox{-2pt}{$\overline{4}$}\\
\cline{1-1}
\raisebox{-2pt}{$\overline{3}$}\\
\cline{1-1}
\raisebox{-2pt}{$\overline{2}$}\\
\cline{1-1}
\end{tabular}~~\longmapsto
\begin{tabular}{|c|c|} \hline
&\raisebox{-2pt}{$4$}\\
\cline{1-2}
 \raisebox{-2pt}{$\overline{4}$}\\
\cline{1-1}
\raisebox{-2pt}{$\overline{3}$}\\
\cline{1-1}
\raisebox{-2pt}{$\overline{2}$}\\
\cline{1-1}
\end{tabular}~~\longmapsto \begin{tabular}{|c|} \hline
 \raisebox{-2pt}{$1$}\\
\hline
\raisebox{-2pt}{$\overline{3}$}\\
\cline{1-1}
\raisebox{-2pt}{$\overline{2}$}\\
\cline{1-1}
\raisebox{-2pt}{$\overline{1}$}\\
\cline{1-1}
\end{tabular}\in QS^{(0,0,0,1)}.$$
Inversement, on a:
$$
 \begin{tabular}{|c|c|c|c|c|} \hline
\raisebox{-2pt}{$~1$}&\raisebox{-2pt}{$~1$}&\raisebox{-2pt}{$~1$}&$\star_5$&$\star_4$\\
\hline
\raisebox{-2pt}{$2$}&\raisebox{-2pt}{$2$}&\raisebox{-2pt}{$\overline{3}$}&$\star_3$&$\star_2$\\
\hline
\raisebox{-2pt}{$3$}&\raisebox{-2pt}{$3$}&\raisebox{-2pt}{$\overline{2}$}&$\star_1$\\
\cline{1-4}
\raisebox{-2pt}{$4$}&\raisebox{-2pt}{$4$}&\raisebox{-2pt}{$\overline{1}$}\\
\cline{1-3}
\end{tabular}~~\stackrel{\sigma}{\longmapsto} \begin{tabular}{cc|c|c|c|} \cline{3-5}
&&\raisebox{-2pt}{$~1$}&\raisebox{-2pt}{$~\overline{4}$}&\raisebox{-2pt}{$~\overline{4}$}\\
\cline{2-5}
& \multicolumn{1}{|c|}{$\star_1$}&\raisebox{-2pt}{$2$}&\raisebox{-2pt}{$\overline{3}$}&\raisebox{-2pt}{$\overline{3}$}\\
\hline
\multicolumn{1}{|c|}{$\star_2$}&$\star_3$&\raisebox{-2pt}{$3$}&\raisebox{-2pt}{$\overline{2}$}&\raisebox{-2pt}{$\overline{2}$}\\
\hline
\multicolumn{1}{|c|}{$\star_4$}&$\star_5$&\raisebox{-2pt}{$\overline{1}$}&\raisebox{-2pt}{$\overline{1}$}&\raisebox{-2pt}{$\overline{1}$}\\
\hline
\end{tabular}~~\longmapsto \begin{tabular}{cc|c|c|c|} \cline{3-5}
&&\raisebox{-2pt}{$~1$}&\raisebox{-2pt}{$~\overline{4}$}&\raisebox{-2pt}{$~\overline{4}$}\\
\cline{2-5}
& \multicolumn{1}{|c|}{}&\raisebox{-2pt}{$2$}&\raisebox{-2pt}{$\overline{3}$}&\raisebox{-2pt}{$\overline{3}$}\\
\hline
\multicolumn{1}{|c|}{\;\;\;}&&\raisebox{-2pt}{$3$}&\raisebox{-2pt}{$\overline{2}$}&\raisebox{-2pt}{$\overline{2}$}\\
\hline
\multicolumn{1}{|c|}{}&$\star_5$&\raisebox{-2pt}{$\overline{1}$}&\raisebox{-2pt}{$\overline{1}$}&\raisebox{-2pt}{$\overline{1}$}\\
\hline
\end{tabular}$$
$$\longmapsto \begin{tabular}{cc|c|c|c|} \cline{3-5}
&&\raisebox{-2pt}{$~2$}&\raisebox{-2pt}{$~\overline{4}$}&\raisebox{-2pt}{$~\overline{4}$}\\
\cline{2-5}
& \multicolumn{1}{|c|}{}&\raisebox{-2pt}{$3$}&\raisebox{-2pt}{$\overline{3}$}&\raisebox{-2pt}{$\overline{3}$}\\
\hline
\multicolumn{1}{|c|}{\;\;\;}&\;\;\;&\raisebox{-2pt}{$4$}&\raisebox{-2pt}{$\overline{2}$}&\raisebox{-2pt}{$\overline{2}$}\\
\hline
\multicolumn{1}{|c|}{$\star_4$}&\raisebox{-2pt}{$\overline{4}$}&\raisebox{-2pt}{$\overline{1}$}&\raisebox{-2pt}{$\overline{1}$}\\
\cline{1-4}
\end{tabular}~~\longmapsto \begin{tabular}{cc|c|c|c|} \cline{3-5}
&&\raisebox{-2pt}{$~2$}&\raisebox{-2pt}{$~\overline{4}$}&\raisebox{-2pt}{$~\overline{4}$}\\
\cline{2-5}
& \multicolumn{1}{|c|}{}&\raisebox{-2pt}{$3$}&\raisebox{-2pt}{$\overline{3}$}&\raisebox{-2pt}{$\overline{3}$}\\
\hline
\multicolumn{1}{|c|}{\;\;\;}&$\star_3$&\raisebox{-2pt}{$4$}&\raisebox{-2pt}{$\overline{2}$}&\raisebox{-2pt}{$\overline{2}$}\\
\hline
\multicolumn{1}{|c|}{\raisebox{-2pt}{$\overline{4}$}}&\raisebox{-2pt}{$\overline{1}$}&\raisebox{-2pt}{$\overline{1}$}\\
\cline{1-3}
\end{tabular}~~\longmapsto \begin{tabular}{cc|c|c|c|} \cline{3-5}
&&\raisebox{-2pt}{$~2$}&\raisebox{-2pt}{$~\overline{4}$}&\raisebox{-2pt}{$~\overline{4}$}\\
\cline{2-5}
& \multicolumn{1}{|c|}{}&\raisebox{-2pt}{$3$}&\raisebox{-2pt}{$\overline{3}$}&\raisebox{-2pt}{$\overline{3}$}\\
\hline
\multicolumn{1}{|c|}{\;\;\;}&\raisebox{-2pt}{$\;4$}&$\star_3$&\raisebox{-2pt}{$\overline{2}$}&\raisebox{-2pt}{$\overline{2}$}\\
\hline
\multicolumn{1}{|c|}{\raisebox{-2pt}{$\overline{4}$}}&\raisebox{-2pt}{$\overline{1}$}&\raisebox{-2pt}{$\overline{1}$}\\
\cline{1-3}
\end{tabular}$$
$$\longmapsto \begin{tabular}{cc|c|c|c|} \cline{3-5}
&&\raisebox{-2pt}{$~2$}&\raisebox{-2pt}{$~\overline{4}$}&\raisebox{-2pt}{$~\overline{4}$}\\
\cline{2-5}
& \multicolumn{1}{|c|}{}&\raisebox{-2pt}{$3$}&\raisebox{-2pt}{$\overline{3}$}&\raisebox{-2pt}{$\overline{3}$}\\
\hline
\multicolumn{1}{|c|}{\;\;\;}&\raisebox{-2pt}{$\;4$}&\raisebox{-2pt}{$\overline{2}$}&$\star_3$&\raisebox{-2pt}{$\overline{2}$}\\
\hline
\multicolumn{1}{|c|}{\raisebox{-2pt}{$\overline{4}$}}&\raisebox{-2pt}{$\overline{1}$}&\raisebox{-2pt}{$\overline{1}$}\\
\cline{1-3}
\end{tabular}~~\longmapsto \begin{tabular}{cc|c|c|c|} \cline{3-5}
&&\raisebox{-2pt}{$~2$}&\raisebox{-2pt}{$~\overline{4}$}&\raisebox{-2pt}{$~\overline{4}$}\\
\cline{2-5}
& \multicolumn{1}{|c|}{}&\raisebox{-2pt}{$3$}&\raisebox{-2pt}{$\overline{3}$}&\raisebox{-2pt}{$\overline{3}$}\\
\hline
\multicolumn{1}{|c|}{$\star_2$}&\raisebox{-2pt}{$\;4$}&\raisebox{-2pt}{$\overline{2}$}&\raisebox{-2pt}{$\overline{2}$}\\
\cline{1-4}
\multicolumn{1}{|c|}{\raisebox{-2pt}{$\overline{4}$}}&\raisebox{-2pt}{$\overline{1}$}&\raisebox{-2pt}{$\overline{1}$}\\
\cline{1-3}
\end{tabular}~~\longmapsto \begin{tabular}{cc|c|c|c|} \cline{3-5}
&&\raisebox{-2pt}{$~2$}&\raisebox{-2pt}{$~\overline{4}$}&\raisebox{-2pt}{$~\overline{4}$}\\
\cline{2-5}
& \multicolumn{1}{|c|}{}&\raisebox{-2pt}{$3$}&\raisebox{-2pt}{$\overline{3}$}&\raisebox{-2pt}{$\overline{3}$}\\
\hline
\multicolumn{1}{|c|}{\raisebox{-2pt}{$\;3$}}&$\star_2$&\raisebox{-2pt}{$\overline{2}$}&\raisebox{-2pt}{$\overline{2}$}\\
\cline{1-4}
\multicolumn{1}{|c|}{\raisebox{-2pt}{$~\overline{3}$}}&\raisebox{-2pt}{$\overline{1}$}&\raisebox{-2pt}{$\overline{1}$}\\
\cline{1-3}
\end{tabular}~$$
$$\longmapsto \begin{tabular}{cc|c|c|c|} \cline{3-5}
&&\raisebox{-2pt}{$~3$}&\raisebox{-2pt}{$~\overline{4}$}&\raisebox{-2pt}{$~\overline{4}$}\\
\cline{2-5}
& \multicolumn{1}{|c|}{}&\raisebox{-2pt}{$4$}&\raisebox{-2pt}{$\overline{3}$}&\raisebox{-2pt}{$\overline{3}$}\\
\hline
\multicolumn{1}{|c|}{$\;3$}&\raisebox{-2pt}{$~\overline{4}$}&$\star_2$&\raisebox{-2pt}{$\overline{2}$}\\
\cline{1-4}
\multicolumn{1}{|c|}{\raisebox{-2pt}{$\overline{3}$}}&\raisebox{-2pt}{$\overline{1}$}&\raisebox{-2pt}{$\overline{1}$}\\
\cline{1-3}
\end{tabular}~~\longmapsto \begin{tabular}{cc|c|c|c|} \cline{3-5}
&&\raisebox{-2pt}{$~3$}&\raisebox{-2pt}{$~\overline{4}$}&\raisebox{-2pt}{$~\overline{4}$}\\
\cline{2-5}
& \multicolumn{1}{|c|}{$\star_1$}&\raisebox{-2pt}{$4$}&\raisebox{-2pt}{$\overline{3}$}&\raisebox{-2pt}{$\overline{3}$}\\
\hline
\multicolumn{1}{|c|}{\raisebox{-2pt}{$\;3$}}&\raisebox{-2pt}{$\overline{4}$}&\raisebox{-2pt}{$\overline{2}$}\\
\cline{1-3}
\multicolumn{1}{|c|}{\raisebox{-2pt}{$\overline{3}$}}&\raisebox{-2pt}{$\overline{1}$}&\raisebox{-2pt}{$\overline{1}$}\\
\cline{1-3}
\end{tabular}~~\longmapsto \begin{tabular}{cc|c|c|c|} \cline{3-5}
&&\raisebox{-2pt}{$~3$}&\raisebox{-2pt}{$~\overline{4}$}&\raisebox{-2pt}{$~\overline{4}$}\\
\cline{2-5}
& \multicolumn{1}{|c|}{$\;3$}&$\star_1$&\raisebox{-2pt}{$\overline{3}$}&\raisebox{-2pt}{$\overline{3}$}\\
\hline
\multicolumn{1}{|c|}{\raisebox{-2pt}{$\;3$}}&\raisebox{-2pt}{$\overline{3}$}&\raisebox{-2pt}{$\overline{2}$}\\
\cline{1-3}
\multicolumn{1}{|c|}{\raisebox{-2pt}{$\overline{3}$}}&\raisebox{-2pt}{$\overline{1}$}&\raisebox{-2pt}{$\overline{1}$}\\
\cline{1-3}
\end{tabular}$$
$$\longmapsto \begin{tabular}{cc|c|c|c|} \cline{3-5}
&&\raisebox{-2pt}{$3$}&\raisebox{-2pt}{$\overline{4}$}&\raisebox{-2pt}{$\overline{4}$}\\
\cline{2-5}
& \multicolumn{1}{|c|}{\raisebox{-2pt}{$3$}}&\raisebox{-2pt}{$\overline{3}$}&\raisebox{-2pt}{$\overline{3}$}\\
\cline{1-4}
\multicolumn{1}{|c|}{\raisebox{-2pt}{$3$}}&\raisebox{-2pt}{$\overline{3}$}&\raisebox{-2pt}{$\overline{2}$}\\
\cline{1-3}
\multicolumn{1}{|c|}{\raisebox{-2pt}{$\overline{3}$}}&\raisebox{-2pt}{$\overline{1}$}&\raisebox{-2pt}{$\overline{1}$}\\
\cline{1-3}
\end{tabular}~~\stackrel{\sigma}{\longmapsto} \begin{tabular}{cc|c|c|c|} \cline{3-5}
&&\raisebox{-2pt}{$1$}&\raisebox{-2pt}{$1$}&\raisebox{-2pt}{$3$}\\
\cline{3-5}
& &\multicolumn{1}{|c|}{$2$}&\raisebox{-2pt}{$3$}&\raisebox{-2pt}{$\overline{3}$}\\
\cline{2-5}
&\multicolumn{1}{|c|}{\raisebox{-2pt}{$3$}}&\raisebox{-2pt}{$3$}&\raisebox{-2pt}{$\overline{3}$}\\
\cline{1-4}
\multicolumn{1}{|c|}{\raisebox{-2pt}{$4$}}&\raisebox{-2pt}{$4$}&\raisebox{-2pt}{$\overline{3}$}\\
\cline{1-3}
\end{tabular}~~\longmapsto \begin{tabular}{|c|c|c|c|c|} \hline
\raisebox{-2pt}{$1$}&\raisebox{-2pt}{$1$}&\raisebox{-2pt}{$1$}&\raisebox{-2pt}{$1$}&\raisebox{-2pt}{$3$}\\
\hline
\raisebox{-2pt}{$2$}&\raisebox{-2pt}{$2$}&\raisebox{-2pt}{$2$}&\raisebox{-2pt}{$3$}&\raisebox{-2pt}{$\overline{3}$}\\
\hline
\raisebox{-2pt}{$3$}&\raisebox{-2pt}{$3$}&\raisebox{-2pt}{$3$}&\raisebox{-2pt}{$\overline{3}$}\\
\cline{1-4}
\raisebox{-2pt}{$4$}&\raisebox{-2pt}{$4$}&\raisebox{-2pt}{$\overline{3}$}\\
\cline{1-3}
\end{tabular}~~\longmapsto \begin{tabular}{|c|c|c|} \hline
\raisebox{-2pt}{$1$} & \raisebox{-2pt}{$1$} & \raisebox{-2pt}{$3$}\\
\hline
\raisebox{-2pt}{$2$} & \raisebox{-2pt}{$3$} & \raisebox{-2pt}{$\overline{3}$}\\
\hline
\raisebox{-2pt}{$3$} & \raisebox{-2pt}{$\overline{3}$}\\
\cline{1-2}
\raisebox{-2pt}{$\overline{3}$}\\
\cline{1-1}
\end{tabular}
.$$
\end{exe}

Comme pour $SL(n)$, on r\'{e}alise le dernier pas de notre preuve en appliquant les
relations de Pl\"ucker internes et externes aux tableaux $T \in NQS^{<\lambda>}$.
On ordonne les tableaux de Young suivant l'ordre habituel : deux tableaux $T$ et $S$ v\'erifient $S<T$ si la forme $\mu=(b_1,\dots,b_n)$ de $S$ est plus petite que la forme $\lambda=(a_1,\dots,a_n)$ de $T$ pour l'ordre lexicographique ou, si ces deux formes sont les m\^emes si, lorsqu'on lit ces deux tableaux colonne par colonne, de droite \`a gauche et dans chaque colonne de bas en haut, le premier couple d'entr\'ees diff\'erentes v\'erifie $t_{i,j}<s_{i,j}$.\\
On fait une r\'ecurrence sur cet ordre total. On suppose que tout tableau $S$ tel que $S<T$ s'\'ecrit modulo les relations de Pl\"ucker et les relations $\delta^{(s)}_{1,2,\dots,s}-1=0$, comme une combinaison lin\'eaire de tableaux quasi standards $U_j$ tels que $U_j\leq S$.\\
Soit maintenant $s$ le plus grand entier tel que $T$ est dans $NQS_s^{<\lambda>}$. Pour tout $\ell$ tel que la colonne $c_\ell$ ait une hauteur $\geq s$, on note $\partial^\ell_j T$ un tableau ayant la m\^eme forme que $T$, dont les colonnes num\'eros $\ell+1,\ell+2,\dots$ sont celles de $T$, le bas de la colonne num\'ero $\ell$ (les entr\'ees des lignes $s+1,\dots$) est le bas de la colonne num\'ero $\ell$ de $T$ et les $s$ premi\`eres entr\'ees de la colonne num\'ero $\ell$ sont $1,2,\dots,s$. Remarquons que $\partial^\ell_j T\geq T$.

On fait une r\'ecurrence sur $\ell$. On suppose que, modulo les relations de Pl\"ucker, il existe des tableaux $\partial^\ell_j T$ et des tableaux $S^\ell_k<T$, de m\^eme forme que $T$, tels que
$$
T=\sum_j \partial^\ell_j T+\sum_{k}S^\ell_k.
$$
Consid\'erons un des tableaux $\partial^\ell_j T$

{\bf{Cas 1:}} $t_{s+1,\ell} > t_{s,\ell+1}$

\

Gr\^ace \`a la relation de Pl\"ucker sur les colonnes num\'eros $\ell$, $\ell+1$, le tableau $\partial^\ell_j T$ s'\'{e}crit:
$$
\partial^\ell_j T=\partial^{\ell+1}_j T+ \displaystyle \sum_{S'<T}~S'.
$$
O\`u $\partial^{\ell+1}_j T$ est obtenu en permutant les $s$ premi\`eres lignes des colonnes num\'eros $\ell$ et $\ell+1$ de $\partial^\ell_j T$.\\

{\bf{Cas 2:}} $t_{s+1,\ell} \leq t_{s,\ell+1}$ et $t_{s+1,\ell}$ est non barr\'{e}.

\

Dans ce cas, puisque $T\in NQS_s^{<\lambda>}$, $t_{s,\ell+1}$ est non barr\'{e}, $\overline{t_{s+1,\ell}}$ appara\^{\i}t
dans la colonne num\'ero $\ell$, on applique une relation interne, sur la colonne num\'ero $\ell$ et sur ce couple d'entr\'{e}es, on obtient :
$$
\partial^\ell_j T=\sum_{i>t_{s,\ell+1}}U_{ij}+\sum_{i<t_{s,\ell+1}}S_{ij},
$$
on remarque que les tableaux $S_{ij}$ sont plus petits que $T$ et d'apr\`es le lemme 3, dans chaque $U_{ij}$, $u_{s+1,\ell}>t_{s,\ell+1}$, en appliquant le cas 1, on peut \'ecrire :
$$
\partial^\ell_j T=\sum_i\partial^{\ell+1}_{i,j} T+ \displaystyle \sum_{S'<T}~S'.
$$

{\bf{Cas 3:}} $t_{s+1,\ell} \leq t_{s,\ell+1}$ et $t_{s+1,\ell}$ est barr\'{e}.

\

Alors l'entr\'{e}e $t_{s,\ell+1}$ est barr\'{e}e, on note $t_{s,\ell+1}=\overline{v}$, alors l'entr\'ee $v$ appara\^\i t dans la colonne $\ell+1$, on applique la relation interne sur la colonne num\'ero $\ell+1$ et sur ce couple d'indice, on obtient :
$$
\partial^\ell_j T=\sum_{i>v}U_{ij}+\sum_{i<v}S_{ij},
$$
on remarque que les tableaux $S_{ij}$ sont plus petits que $T$ et d'apr\`es le lemme 3, dans chaque $U_{ij}$, $u_{s,\ell+1}<t_{s+1,\ell}$, de plus $u_{s',\ell+1}=t_{s',\ell+1}$, pour tout $s'>s$. En appliquant le cas 1, on peut \'ecrire :
$$
\partial^\ell_j T=\sum_i\partial^{\ell+1}_{i,j} T+ \displaystyle \sum_{S'<T}~S'.
$$

Donc par r\'ecurrence, on a bien, pour tout $\ell$ tel que la colonne $c_\ell$ de $T$ ait au moins $s$ cases,
$$
T=\sum_j \partial^\ell_j T+\sum_k S_k,
$$
avec $S_k<T$ pour tout $k$. On \'ecrit cette relation pour la premi\`ere colonne de hauteur $s$ de $T$, on obtient des tableaux $\partial^\ell_j T$ ayant une colonne triviale, on supprime cette colonne triviale gr\^ace \`a la relation $\delta^{(s)}_{1,2,\dots,s}=1$, on obtient des tableaux $(\partial^\ell_j T)'<T$. Donc $T$ est une combinaison lin\'eaire de tableaux $S<T$, par induction, c'est une combinaison lin\'eaire de tableaux quasi standards de forme $\mu\subset\lambda$. 

L'ensemble $\cup _{\mu \subset \lambda}~~QS^{<\mu>}$ est un syst\`eme g\'en\'erateur du $N^+$ module $\mathbb S^{<\lambda>}_{\mathfrak{n}_{\mathfrak{sp}(2n)}^+}\subset \mathbb S^{<\lambda>}_{red}$, ce module a pour dimension le cardinal de $SS^{<\lambda>}$, le syst\`eme $\cup _{\mu \subset \lambda}~~QS^{<\mu>}$ est donc une base de ce module.

\begin{theo}

\

Tout tableau de $SS^{<\lambda>}$ est une combinaison lin\'{e}aire de
tableaux de $\cup _{\mu \subset
\lambda}~~QS^{<\mu>}$.\\
L'ensemble $QS^{<\bullet>}$ est une base de
$\mathbb{S}^{<\bullet>}$, adapt\'{e}e \`{a} la stratification des
$N^+$-modules $\mathbb{S}^{<\lambda>}_{\mathfrak{n}_{\mathfrak{sp}(2n)}^+}$.

\end{theo}



\begin{thebibliography}{ZZZZZ}
\bibitem [AAK]{AAK}
B. Agrebaoui, D. Arnal, O. Khlifi, ``Diamant representations of rank
two semisimple Lie algebras''; A para\^\i tre, Journal of Lie theory (2008).

\bibitem[ABW] {ABW} D. Arnal, N. Bel Baraka, N. Wildberger :
``Diamond representations of $\mathfrak{sl}(n)$'', Ann. Math. Blaise
Pascal, 13 $\hbox{n} ^{\circ } 2$ (2006), p.381--429.

\bibitem [D]{D} R. G. Donnelly, ''Explicit Constructions of the fundamental reprentations
 of the symplectic Lie algebras''; Journal of algebra 223, p.37-64, (2000).

\bibitem [DeC]{DeC} C. De Concini, ''Symplectic standard tableaux'',
Advances in Math. 34 (1979), p.1-27, MR80m:14036.

\bibitem [FH]{FH}
W. Fulton and J. Harris, ``Representation theory''; Readings in
Mathematics 129 (1991) Springer- Verlag, New York.

\bibitem [H]{H} J. E. W. Humphreys, ``Introduction to Lie algebras and
representation theory''; Springer-Verlag, New York; Heidelberg;
Berlin (1972).

\bibitem [KN]{KN} M. Kashiwara, T. Nakashima, ''Crystal graphs for
representations of the q-analogue of classical Lie algebras'',
Journal of algebra 165 (1994), p.295-345.

\bibitem [L]{L} C. Lecouvey, ''Kostka-Foulkes polynomials cyclage graphs and charge statistic for
the root system $C_n$''; Journal of Algebraic Combinatorics 21, pp.
203-240 (2005).

\bibitem [SH]{SH} J. T. Sheats, '' A symplectic jeu de taquin
bijection between the tableaux of King and of De Concini'';
Transaction of the American Mathematical Society, volume 351, Number
9, p.3569-3607, S 0002-9947(99)02166-2 (1999).

\bibitem [V]{V}  V.S. Varadarajan  ``Lie groups, Lie algebras, and their
representations``; Springer-Verlag, New York ; Berlin (1984).

\bibitem [W]{W}  N.J. Wildberger  ``Quarks, diamond and representations of $\mathfrak{sl}(3)$``; Preprint.




\end{thebibliography}
\end{document}